\documentclass[,11pt,leqno]{article}
\usepackage{latexsym}
\usepackage{amsfonts}
\usepackage{amsmath}
\newcommand\blackslug{\hbox{\hskip 1pt \vrule width 4pt height 8pt depth 1.5pt
        \hskip 1pt}}
\newcommand\bbox{\hfill \quad \blackslug \bigbreak}

\def\ll{,\ldots,}

\makeatletter
\newcommand{\leqnomode}{\tagsleft@true}
\newcommand{\reqnomode}{\tagsleft@false}
\makeatother

\newcommand{\sset}[1]{\left\{#1\right\}}
\newcommand{\vare}{\varepsilon}
\renewcommand{\epsilon}{\varepsilon}

\newcommand{\floor}[1]{\left\lfloor#1\right\rfloor}
\newcommand{\ceil}[1]{\left\lceil#1\right\rceil}

\title{Pure pairs. III. Sparse graphs with no polynomial-sized anticomplete pairs}
\author{Maria Chudnovsky\thanks{Supported by NSF grant DMS-1550991.
This material is based upon work supported in part by the U. S. Army
Research Laboratory and the U. S. Army Research Office under grant
number
W911NF-16-1-0404.}\\
Princeton University, Princeton, NJ 08544
\\
\\
Jacob Fox\\
Stanford University, Stanford, CA 94305-2125
\\
\\
Alex Scott\thanks{Supported by a Leverhulme Trust Research
Fellowship}\\
Mathematical Institute, University of Oxford, Oxford OX2 6GG, UK
\\
\\
Paul Seymour\thanks{Supported by ONR grant N00014-14-1-0084 and NSF
grant DMS-1265563.}\\
Princeton University, Princeton, NJ 08544
\\
\\
Sophie Spirkl\\
Princeton University, Princeton, NJ 08544}

\date{December 31, 2017; revised \today}

\newtheorem{thm}{}[section]

\newcommand{\Proof}{\noindent{\bf Proof.}\ \ }

\begin{document}
\maketitle
\begin{abstract} A graph is {\em $H$-free} if it has no induced subgraph isomorphic to $H$, and $|G|$ denotes the number of vertices of $G$.
A conjecture of Conlon, Sudakov and the second author asserts that:
\begin{itemize}
\item For every graph $H$, there exists $\vare>0$
such that in every $H$-free graph $G$ with $|G|>1$ there are
two disjoint sets of vertices, of sizes at least $\vare |G|^\vare$ and $\vare |G|$, complete or anticomplete
to each other. 
\end{itemize}
This is equivalent to:
\begin{itemize}
\item The ``sparse linear conjecture'': For every graph $H$, there exists $\vare>0$
such that in every $H$-free graph $G$ with $|G|>1$, either some vertex
has degree at least $\vare |G|$, or there are two disjoint sets of vertices, of sizes
at least $\vare |G|^\vare$ and $\vare |G|$, anticomplete to each other. 
\end{itemize}
We prove a number of partial results towards the sparse linear conjecture. In particular,
we prove it holds for a large class of graphs $H$, and we prove that something like it holds for all graphs $H$. 
More exactly, say $H$ is ``almost-bipartite'' if $H$ is triangle-free and $V(H)$ can be partitioned 
into a stable set and a set inducing a graph of maximum degree at most one. (This includes 
all graphs that arise from another graph by subdividing every edge at least once.) Our main result is:

\begin{itemize}
\item The sparse linear conjecture holds
for all almost-bipartite graphs $H$. 
\end{itemize}
(It remains open when $H$ is the triangle $K_3$.) There is also a stronger theorem:
\begin{itemize}
\item For every almost-bipartite graph $H$, there exist $\vare, t >0$
such that for every graph $G$ with $|G|> 1$ and maximum degree less than $\vare |G|$, and for every $c$ with $0<c\le 1$, either
$G$ contains $\vare c^t |G|^{|H|}$  induced copies of $H$, or there are two disjoint sets $A, B \subseteq V(G)$ 
with $|A| \geq \vare c^t |G|$ and $|B| \geq \vare |G|$, and with at most $c|A|\cdot |B|$ edges between them.
\end{itemize}
We also prove some variations on the sparse linear conjecture, such as:
\begin{itemize}
\item For every graph $H$, there exists $\vare >0$
such that in every $H$-free graph $G$ with $|G|>1$ vertices, either some vertex
has degree at least $\vare |G|$, or there are two disjoint sets $A, B$ of vertices with $|A|\cdot |B| \geq \vare |G|^{1 + \vare}$,
anticomplete to each other. 
\end{itemize}

\end{abstract}

\section{Introduction}

All graphs in this paper are finite and have no loops or parallel edges.
We denote the number of vertices of $G$ by $|G|$, and $G[X]$ denotes the
subgraph induced on $X\subseteq V(G)$.
If $G,H$ are graphs,
we say that $G$ {\em contains} $H$ if some induced subgraph of $G$ is isomorphic to $H$, and $G$ is
{\em $H$-free} otherwise. 
If $A,B\subseteq V(G)$ are disjoint, we say $A$ is {\em complete} to $B$ if every vertex in $A$ is adjacent to every
vertex in $B$, and {\em anticomplete} to $B$ if there is no edge between $A$ and $B$. A pair $(A,B)$ of subsets of $V(G)$ is {\em pure}
if $A\cap B=\emptyset$ and $A$ is either complete or anticomplete to $B$.

Erd\H{o}s, Hajnal and Pach~\cite{EHP} proved:

\begin{thm}\label{EHP}
For every graph $H$ there exists $\vare >0$  such that for every $H$-free
graph $G$ with $|G|>1$, there is a pure pair $(A,B)$, with 
$|A|,|B|\ge \vare |G|^\vare$. 
\end{thm}

The goal of this paper is to strengthen \ref{EHP}. For instance, we shall
prove in section \ref{sec:ehpimprove} that:

\begin{thm}\label{easyEH}
For every graph $H$ there exists $\vare >0$ such that in every $H$-free graph $G$ with $|G|>1$, 
there is a pure pair $(A,B)$ 
with $|A|\cdot|B|\ge \vare |G|^{1+\vare}$.
\end{thm}

This is best possible up to $\vare$, as can be seen from an Erd\H{o}s-Renyi
random graph $G(n,p)$ with $n$ very large, $p=n^{-\delta}$ and $\delta<\vare$, taking $H$ to be a complete graph
on more than $1+2/\delta$ vertices.
A union bound shows that, almost surely, this random graph is $H$-free,
and does not contain any such pair $(A,B)$.

We will also prove a number of other strengthenings of \ref{EHP}. Here are some ways in which we could try to modify it:
\begin{itemize}
\item Make $|B|\ge \vare |G|$; conjecture \ref{foxsymm} below says this is always possible. We cannot ask for both $A,B$ to be linear, however;
a random graph construction shows that
this can only be true when both $H$ and its complement $\overline{H}$ are forests, that is, $H$
is an induced subgraph of a four-vertex path.
\item Get more than two sets.
\item Replace ``$H$-free'' by a weaker hypothesis, that there are not many copies of $H$ in $G$.
\item Generalize ``pure'' to ``$(1-c)$-dense or $c$-sparse'' (these limit the number of edges between $A,B$).
\item Assuming that $G$ is ``$\vare$-bounded'' (that is, its maximum degree is less than $\vare |G|-1$), eliminate the ``complete'' or ``$(1-c)$-dense''
outcome. 
\end{itemize}
We will prove various combinations of these. For instance, in \ref{EHdensityimprove} we satisfy the second, third and fourth bullets, and also
the fifth in \ref{lem:splusone}. Our main result, \ref{almostbip}, satisfies the first, third, fourth and fifth bullets, but only when $H$
is ``almost-bipartite''.

We need a few definitions. Let $G$ be a graph. 
For a vertex $v \in V(G)$, we use $N(v)$ to denote its set of neighbours, and we define $N[v] = N(v) \cup \sset{v}$.
For $\vare>0$, a graph $G$ is \emph{$\vare$-bounded} if $|N[v]| < \vare |G|$ for all $v \in V(G)$. 
A pair $(A, B)$ of subsets of $V(G)$ is
an \emph{$(x, y)$-pair} for $x,y\ge 0$, if $A\cap B=\emptyset$ and $|A| \geq x$ and $|B| \geq y$; and if $(A,B)$
is also pure we call it a {\em pure \emph{$(x, y)$-pair}}.

There is a conjecture of Conlon, Sudakov and the second author (problem 3.13 in~\cite{fox}) that more than \ref{easyEH} is true,
that we can make the larger of $A,B$ linear in $n$:

\begin{thm}\label{foxsymm}
{\bf Conjecture: }For every graph $H$ there exists $\epsilon > 0$ such that every $H$-free graph $G$ 
with $|G|>1$ contains a pure $(\vare |G|^\vare, \vare |G|)$-pair. 
\end{thm}

If we restrict attention to $\vare$-bounded graphs,
then there cannot exist a pair
as in \ref{foxsymm} complete to each other; so 
\ref{foxsymm} would imply:

\begin{thm}\label{ourconj}
{\bf Conjecture: }For every graph $H$ there exists $\epsilon > 0$ such that in every $H$-free $\vare$-bounded graph $G$
with $|G|> 1$ vertices, there is an anticomplete $(\vare |G|^\vare, \vare |G|)$-pair.
\end{thm}

A theorem of R\"odl~\cite{rodl}
shows that a graph $H$ satisfies \ref{foxsymm} if and only both $H$ and $\overline{H}$
satisfy \ref{ourconj}. Thus \ref{foxsymm} (for all $H$) is equivalent to \ref{ourconj} (for all $H$). On the other hand, 
for certain graphs $H$, \ref{ourconj} turns out to be much more tractable than \ref{foxsymm}. 

We say $H$ is      
{\em almost-bipartite} if it is triangle-free and its vertex set can be partitioned into a stable set and a 
set that induces a subgraph with maximum degree at most one.
A consequence of our main result is:
\begin{thm}\label{almostbip0}
All almost-bipartite graphs $H$ satisfy \ref{ourconj}.
\end{thm}
(It remains open whether $H=K_3$ satisfies \ref{ourconj}, however.)
The full conjecture \ref{foxsymm} has not been proved for many graphs $H$.
In \cite{ehfive}, the authors prove that \ref{foxsymm} holds for a five-cycle, but otherwise it has only been proved 
so far for graphs $H$ that are induced subgraphs of a four-vertex path.
A consequence of our results is that 
two more graphs satisfy
\ref{foxsymm}, namely a four-cycle and its complement.

\section{A digression}

In order to further motivate \ref{foxsymm}, let us mention in passing that a stronger,  related conjecture  is false. It would have been nice if:

\begin{thm}\label{falseconj}
{\bf False conjecture: } For every graph $H$ there exist $c>0$ and a function $f$ with $f(n)=o(n)$, such that for every $H$-free graph $G$ with $|G|>1$, there are disjoint
subsets $A,B$ of $V(G)$ with $|A|, |B|\ge c|G|$, such that either
\begin{itemize}
\item every vertex in $A$ has at most $f(|G|)$ neighbours in $B$ and vice versa, or
\item every vertex in $A$ has at most $f(|G|)$ non-neighbours in $B$ and vice versa.
\end{itemize}
\end{thm}
This would have been nice for two reasons, at least if $f(n)$ could be taken to be $n^{1-\vare}$ for some $\vare>0$, because then
it implies the Erd\H{o}s-Hajnal conjecture~\cite{EH0,EH}, and
it implies \ref{foxsymm}. 
But it is false, for all $f$: here is a sketch of a counterexample. Let $H$ be a cycle of length $g\ne 4$,
let $c>0$, let $f$ be a function with $f(n)=o(n)$, and let us show that
$c,f$ do not satisfy \ref{falseconj}.
Choose $d$ with $(g-1)/g<d<1$.
Let $m$ be a large integer.
Let $F$ be an $H$-free graph on $m$ vertices with the property that, for every two disjoint subsets $A,B$ of $V(F)$ of cardinality at
least $cm/4$,  the number of edges between $A,B$
is at least $\frac12 m^{-d}|A|\cdot |B|$ and at most $\frac12 |A|\cdot|B| $. (If $m$ is large enough, such a graph $F$ exists, since with high probability
the Erd\H{o}s-R\'enyi random graph $G(2m,p)$ with $p=m^{-d}$ has an induced subgraph on $m$ vertices with this property.)

Now let $n$ be a large multiple of $m$. Blow-up
$F$ to obtain $G$ by replacing each vertex $i$ by a stable set $V_i$ of cardinality
$n/m$, and make $V_i$ complete to $V_j$ if $i$ is adjacent to $j$, and anticomplete otherwise.
It is fairly easy to check that $G$ is $H$-free (this is where we use that $g\ne 4$), and for all disjoint subsets $A,B$ of $V(G)$ both of
cardinality at least $cn$, the
number of edges of $G$ between $A,B$ is at least $\frac{c^2}{128}m^{-d}|A|\cdot|B|$
and the number of nonedges between $A,B$ is at least $\frac{1}{128}|A|\cdot|B|$; and both of these are more than $f(n)n$ if $n/m$ is  
large enough. Thus $G$ is a counterexample.

To make a counterexample when $H$ is a cycle of length four, we make the same construction, except we take the $V_i$'s to be cliques. Thus
\ref{falseconj} is false for every $H$ that has an induced cycle, and therefore for every $H$ whose complement has an induced cycle.

\section{Density theorems} \label{sec:density}

If $A,B\subseteq V(G)$ are disjoint, $E(A,B)$ denotes the set of edges of $G$
with one end in $A$ and one in $B$.
For $c\ge 0$, a pair $(A, B)$ of subsets of $V(G)$ is
\begin{itemize}
\item \emph{$c$-sparse} if $A\cap B=\emptyset$ and $|E(A, B)| \leq c|A|\cdot |B|$; and
\item \emph{$c$-dense} if $A\cap B=\emptyset$ and $|E(A,B)| \geq c |A|\cdot |B|$.
\end{itemize}

The best general bound for the Erd\H{o}s-Hajnal conjecture~\cite{EH0} to date was proved by Erd\H{o}s and Hajnal in~\cite{EH},
namely:

\begin{thm}\label{EHpartial}
For every graph $H$, there exists $\vare>0$ such that for every non-null $H$-free graph $G$,
some clique or stable set of $G$ has cardinality at least $2^{\vare \sqrt{\log |G|}}$.
\end{thm}

One of the key steps in proving this was to prove the following:
\begin{thm}\label{EHdensity}
For every graph $H$, there exist $\vare,s>0$ such that for every $H$-free graph $G$ with $|G|>1$, 
and every $c$ with $0\le c\le 1$, $G$ contains either a $c$-sparse or a $(1-c)$-dense $(\vare c^s |G|, \vare c^s |G|)$-pair. 
\end{thm}

Conlon, Sudakov and the second author (problem 3.13 in ~\cite{fox}) asked whether one of the sets $A,B$ could always be chosen of linear size, 
independent of $c$: that is,

\begin{thm}\label{foxdensity}
{\bf Conjecture: }For every graph $H$ there exist $\vare,s > 0$ such that
for every $H$-free graph $G$ with $|G|>1$, and all $c$ with $0\le c\le 1$, $G$ contains 
either a $c$-sparse or a $(1-c)$-dense $(\vare c^s |G|, \vare |G|)$-pair.
\end{thm}
For $\vare$-bounded graphs this becomes:

\begin{thm}\label{foxweak}
{\bf Conjecture: }For every graph $H$ there exist $\vare, s > 0$ such that
for every $H$-free $\vare$-bounded graph $G$ with $|G|>1$, and all $c$ with $0\le c\le 1$, there is a 
$c$-sparse $(\vare c^s |G|, \vare |G|)$-pair in $G$. 
\end{thm}

There are implications between these conjectures; in fact
\ref{foxdensity} $\Rightarrow$ \ref{foxsymm} $\Rightarrow$ \ref{ourconj} and
\ref{foxdensity} $\Rightarrow$ \ref{foxweak} $\Rightarrow$ \ref{ourconj}.
We have seen that \ref{foxsymm} implies \ref{ourconj}.

\bigskip

\noindent{\bf Proof of \ref{foxweak}, assuming \ref{foxdensity}.\ \ }
Let 
$\epsilon'$ and $s$ satisfy \ref{foxdensity}, and let $\epsilon = \epsilon'/2$.
Let $G$ be $H$-free,
and let $0<c\le 1$, and $n=|G|\ge 2$.
By \ref{foxdensity}, there exist disjoint $A,B$ with
$|A|\ge \vare'c^{s} n$
and $|B|\ge \epsilon' n$  such that $(A, B)$ is either $c$-sparse or $(1-c)$-dense. If $(A,B)$ is a $c$-sparse pair, then $(A,B)$ satisfies
\ref{foxweak} as required, because 
$\epsilon' \ge \epsilon $.
If $(A, B)$ is $(1-c)$-dense and not $c$-sparse, and so $\max(c,1-c)$-dense and therefore at least $(1/2)$-dense, 
then some vertex in $A$ has degree at least
$|B|/2\ge \epsilon n$ and again \ref{foxweak} holds.~\bbox

\noindent{\bf Proof of \ref{ourconj}, assuming \ref{foxweak}.\ \ }
Let $\epsilon'$ and $s$ satisfy \ref{foxweak}; and 
let 
$$\vare= \min\left(\epsilon'/2, 1/(s+1),1/4\right).$$ 
We claim that $\epsilon$ satisfies \ref{ourconj}.

Let
$G$ be $\vare$-bounded and $H$-free with $n>1$ vertices, and let $x= \epsilon n^{\epsilon}$.
Choose $c$ such that $c^{s} n = n^{\epsilon}$, that is,
$c = n^{-(1-\epsilon)/s}$. 
\\
\\
(1) {\em We may assume that $x\ge 1$, and $x\le \vare'c^{s} n$, and $cx\le 1/4$.}
\\
\\
Let $v\in V(G)$.
Since $|N[v]|<\epsilon n$, and $\epsilon\le 1/2$, it follows that $v$ has at least $\vare n$ non-neighbours;
and since we may assume that $v$ and its non-neighbours do not form an anticomplete $(\vare n^{\vare}, \vare n)$-pair,
it follows that $x\ge 1$. The second claim holds since $c^{s} n = n^{\vare}$; and the third since 
$cx=n^{-(1-\epsilon)/s}(\epsilon n^{\epsilon})$, and $\vare -(1-\vare)/s\le 0$ and $\vare\le 1/4$.
This proves (1).

\bigskip

Now $c\le 1$, so 
by \ref{foxweak},
there is a $c$-sparse $(\vare'c^{s} n, \vare' n)$-pair $(A, B)$.
By (1) and \ref{sparsetoanti} below, there is an anticomplete $(x,|B|/2)$-pair $(A', B')$ with $A'\subseteq A$ and $B'\subseteq B$.
But then $(A', B')$ satisfies \ref{ourconj}.
This proves \ref{ourconj}.~\bbox

(The proof that \ref{foxdensity} implies \ref{foxsymm} is similar and we omit it.)
We just used a lemma that produces anticomplete pairs from $c$-sparse pairs:
\begin{thm}\label{sparsetoanti}
Let $(A,B)$ be a $c$-sparse pair in a graph $G$. If $x\ge 1/2$ (not necessarily an integer), and 
$x\le |A|$ and $cx\le 1/4$, there is an anticomplete $(x,|B|/2)$-pair $(A',B')$ with $A'\subseteq A$ and $B'\subseteq B$.
\end{thm}
\Proof
Let $d=\ceil{x}$; then since $x\ge 1/2$ it follows that $d\le 2x$. Since $x\le |A|$ and hence $d\le |A|$, there is a subset of $A$ 
with cardinality $d$. By averaging over all
such subsets, it follows that there exists $A'\subseteq A$ with $|A'|=d$ such that $(A', B)$ is $c$-sparse. In particular, there are
at most $cd|B|\le 2cx|B|\le |B|/2$ vertices in $B$ with a neighbour in $A'$; let $B'$ be the other vertices in $B$, and then the result holds.~\bbox

\section{Saturation} \label{sec:counting}

For a graph $H$ and a graph $G$, a \emph{copy of $H$ in $G$} is an isomorphism $\phi$ 
between $H$ and an induced subgraph of $G$. (Thus, there are six copies of $K_3$ in $K_3$.)
In particular, $G$ contains $H$ 
if and only if there is a copy of $H$ in $G$. For $\alpha \geq 0$, we say that a graph $G$ is \emph{$(\alpha, H)$-saturated} if 
there are at least $\alpha |G|^{|H|}$ copies of $H$ in $G$.

With all these results and conjectures, one can try replacing ``$H$-free'' by ``not $(\alpha, H)$-saturated'' for the appropriate choice of $\alpha$.
For instance, we mentioned earlier a theorem of R\"odl~\cite{rodl}; it says that for all $H$ and $\vare>0$, there exists $\delta>0$
such that if $G$ is $H$-free, there is an induced subgraph $J$ with $|J|\ge \delta|G|$ such that 
one of $|E(J)|,|E(\overline{J})|$ is at most $\vare |J|(|J| -1)/2$.
There is a saturation version of this,
the following, as was observed by Sudakov and the second author~\cite{foxsudakov}: 
\begin{thm}\label{thm:copies}
Let $H$ be a graph, and let $\vare > 0$. Then there exist $\alpha, \delta > 0$ such that for every graph $G$, 
if $G$ is not $(\alpha, H)$-saturated, then $G$ contains an induced subgraph $J$ with $|J| \geq \delta |G|$  such that 
one of $|E(J)|,|E(\overline{J})|$ is at most $\vare |J|(|J| -1)/2$. 
\end{thm}

Similarly, one can strengthen \ref{EHdensity}. In fact we will prove the following in section \ref{sec:ehpimprove}; it strengthens \ref{EHdensity} in two ways,
replacing ``$H$-free'' by ``not $(\alpha, H)$-saturated'' and producing $k$ sets instead of two.
($\mathbb{N}$ denotes the set of non-negative integers.)
\begin{thm}\label{EHdensityimprove}
For every graph $H$ and $k \in \mathbb{N}$, there exist $\vare, s, K >0$ such that for every graph $G$ with $n>K$ vertices, 
and every $c$ with $0\le c\le 1$, if $G$ is not $(\vare c^s, H)$-saturated, then there are pairwise disjoint subsets $A_1, \dots, A_k\subseteq V(G)$ such that either:
\begin{itemize}
\item $(A_i, A_j)$ 
is a $c$-sparse $(\vare c^s n, \vare c^s n)$-pair for all distinct $i, j \in \sset{1, \dots, k}$; or
\item $(A_i, A_j)$ 
is a $(1-c)$-dense $(\vare c^s n, \vare c^s n)$-pair for all distinct $i, j \in \sset{1, \dots, k}$.
\end{itemize}
\end{thm}
As usual, as we will prove in \ref{lem:splusone}, if we require that $G$ is $\vare$-bounded, then we can omit the second outcome. 

In light of this, one might try the saturation strengthenings of the two conjectures from the previous section.
\ref{foxdensity} could be strengthened to:
\begin{thm}\label{foxdensitysat}
{\bf Conjecture: }For every graph $H$ there exist $\epsilon,s > 0$ such that
for every graph $G$ with $|G|>1$, and all $c$ with $0\le c\le 1$, either $G$ is $(\vare c^s, H)$-saturated, or 
there is a $c$-sparse or a $(1-c)$-dense $(\vare c^s |G|, \vare |G|)$-pair in $G$.
\end{thm}
Similarly, \ref{foxweak} could be strengthened to:

\begin{thm}\label{foxweaksat}
{\bf Conjecture: }For every graph $H$ there exist $\epsilon, s > 0$ such that for every $\vare$-bounded graph $G$ and all $c$ with $0\le c\le 1$, 
either $G$ is $(\vare c^s, H)$-saturated, or there is a 
$c$-sparse $(\vare c^s |G|, \vare |G|)$-pair in $G$. 
\end{thm}

As in section \ref{sec:density}, we have the implication \ref{foxdensitysat} $\Rightarrow$ \ref{foxweaksat}, and clearly
\ref{foxdensitysat} $\Rightarrow$ \ref{foxdensity}, and \ref{foxweaksat} $\Rightarrow$ \ref{foxweak}. Moreover,
\ref{thm:copies} shows that $H$ satisfies \ref{foxdensitysat} if and only if both $H$ and $\overline{H}$ satisfy \ref{foxweaksat}.

We will prove \ref{foxweaksat} when $H$ is almost-bipartite. Our main theorem (proved in section \ref{sec:almostbip})
says:
\begin{thm}\label{mainthm}
For every almost-bipartite graph $H$ there exist $\epsilon, s > 0$ such that for every $\vare$-bounded graph $G$
and all $c$ with $0\le c\le 1$, either $G$ is $(\vare c^s, H)$-saturated, or there is a
$c$-sparse $(\vare c^s |G|, \vare |G|)$-pair in $G$.
\end{thm}

The graphs $H$ such that both $H$ and $\overline{H}$ are almost-bipartite are the five-cycle, the four-cycle, and its complement, 
as well as all induced subgraphs of these graphs. Therefore, our results imply that \ref{foxdensitysat} holds for these graphs. 

Finally, we remark that we cannot do better than ``$(\vare c^s, H)$-saturated'' in \ref{foxweaksat}, that is, \ref{foxweaksat} 
becomes false if we replace ``$(\vare c^s, H)$-saturated'' by ``$(\vare, H)$-saturated''. This can be seen by letting $H = K_2$.
Let $\vare, s>0$; we will show they do not satisfy the modified \ref{foxweaksat}.
Let $n \in \mathbb{N}$, $\delta = 1/(2s+2)$, and $p = n^{-\delta/2}$,
and let $G$ be an $n$-vertex random graph in which every edge is present independently with 
probability $p$. It follows that $G$ has $\approx \frac12 n^{2-\delta/2}$ edges in expectation, so for $n$ sufficiently large, 
with high probability $G$ 
is not $(\vare, H)$-saturated. Also the probability that there is an anticomplete $(n^\delta, \frac12 \vare n)$-pair in $G$ 
is at most
$$3^n(1-p)^{\frac12 \vare n^{1+\delta}} \leq 3^n e^{-\frac12 \vare n^{1+\delta/2}} \rightarrow 0$$
as $n \rightarrow \infty$; so for $n$ large, with high probability, $G$ has no anticomplete $(n^\delta, \frac12 \vare n)$-pair.

Let $c = n^{-\delta}/4$.
Since $n^\delta\ge 1$, and $n^\delta\le \vare c^s n$ (for large $n$), and $cn^{\delta}= 1/4$, it follows from 
\ref{sparsetoanti} (with $x=n^\delta$) that, if there is no anticomplete $(n^\delta, \frac12 \vare n)$-pair in $G$, then
 there is also no $c$-sparse $(\vare c^s n, \vare n)$-pair in $G$. 
So with high probability, $G$ is not $(\vare, H)$-saturated and has no
$c$-sparse $(\vare c^k n, \vare n)$-pair.

\section{A game on a graph} \label{sec:triple}

To prove our main results, we will use the existence of a winning strategy in a certain game, that we discuss in this section.
Let $G$ be a graph, and let $T \subseteq V(G)$ be a stable set. A graph $G'$ is a \emph{$T$-successor of $G$} if
\begin{itemize}
\item $V(G) = V(G')$ and $G$ is a proper subgraph of $G'$; and
\item every edge in $E(G') \setminus E(G)$ has both ends in $T$.
\end{itemize}

Let $H$ be a graph. For $k \geq 2$ and $m \geq k$, the \emph{$k$-tuple game for $H$} on $m$ vertices is the following game between two players, A and B. 
Let $G_0$ be a graph with $m$ vertices and no edges. Rounds of the game will add edges to $G_0$, making a sequence of graphs $G_1,G_2,\cdots$, all with the same
vertex set and each a proper subgraph of the next. In round $i$, player A selects a stable set $T$ of cardinality $k$ in $G_{i-1}$, and player B
choose a $T$-successor $G_{i}$ of $G_{i-1}$.
Player A wins if at some stage there is an induced subgraph isomorphic 
to $H$.

More precisely, the $i$th round (starting with $i=1$) consists of the following:
\begin{itemize}
\item if $G_{i-1}$ contains $H$, player A has won; otherwise, if $G_{i-1}$ has no stable set of cardinality $k$,
then player B has won;
\item if neither of these, player A chooses a $k$-vertex stable subset $T$ of $G_{i-1}$, and
player B chooses a $T$-successor $G_i$ of $G_{i-1}$.
\end{itemize}
Then a new round commences. Since at least one edge is added in every round, this game terminates after a finite number of rounds.

For a graph $H$, we say that $H$ is \emph{$(m, k)$-forcible} if there is a strategy for player A to play the $k$-tuple game for $H$
on $m$ vertices and always win, that is, reach a graph $G$ that contains $H$. We say that such a strategy \emph{forces} $H$. 
The main result of this section is that for every $H$ and $k$, there exists $m\ge 0$ such that $H$ is $(m, k)$-forcible. 
We begin by proving the base cases: 

\begin{thm} \label{lem:basetriple}
Every graph $H$ is $(|H|,2)$-forcible. If furthermore $|E(H)| = 0$, then $H$ is $(|H|, k)$-forcible 
for all $k \geq 2$. 
\end{thm}
\Proof
The second statement holds since $G_0$ is isomorphic to $H$ if $H$ has no edges. For the first statement, player A picks a 
bijection $f$ between $V(H)$ and $V(G_0)$, and player A ensures that in every round $i$, $G_{i-1}$ is isomorphic to a 
(not necessarily induced) subgraph of $H$. Therefore either $H$ is isomorphic to $G_{i-1}$ or there is an edge $uv \in E(H)$ 
such that $f(u),f(v)$ are not adjacent in $G_{i-1}$. In the first case, player A wins the game; in the second, player A 
picks the set $\sset{f(u), f(v)}$. This forces player B to add the edge $f(u)f(v)$. After $|E(H)|$ rounds, $G_i$ 
has $|E(H)|$ edges and hence is isomorphic to $H$.~\bbox

\begin{thm} \label{lem:union}
Let $H_1, H_2$ be graphs, and let $k \geq 2$ and $m_1, m_2 \in \mathbb{N}$ such that for all $i \in \sset{1,2}$, $H_i$ is 
$(m_i, k)$-forcible. Then the disjoint union of $H_1$ and $H_2$ is $(m_1 + m_2, k)$-forcible.
\end{thm}
\Proof
Let $G_0$ have $m_1 + m_2$ vertices, partitioned into $V_1, V_2$ with $|V_i| = m_i$ for $i \in \sset{1, 2}$. 
Player A first plays according to the $k$-tuple game for $H_1$ on $G_0[V_1]$. Since $H_1$ is $(m_1, k)$-forcible, this game 
stops after a finite number $s_1$ of rounds and $G_{s_1-1}[V_1]$ contains $H_1$. Since every $k$-tuple 
picked by player A is contained in $V_1$, it follows that $V_2$ is stable and anticomplete to $V_1$ in $G_{s_1-1}$. 
Instead of stopping in round $s_1$, player A now plays according to the $k$-tuple game for $H_2$ on $G_{s_1-1}[V_2]$. 
After a finite number $s_2$ of rounds, $G_{s_1+s_2-1}[V_2]$ contains $H_2$, and $V_1$ remains 
anticomplete to $V_2$. But this implies that $G_{s_1 + s_2 - 1}$ contains the 
disjoint union of $H_1$ and $H_2$.~\bbox

\begin{thm} \label{lem:triplemain}
Let $H$ be a graph and $k \geq 2$. Then there exists $m \geq k$ such that $H$ is $(m, k)$-forcible. 
\end{thm}

\Proof
We prove this by induction on $k+|E(H)|$. The statement holds in the base cases, when either $k = 2$ or $|E(H)| = 0$, 
by \ref{lem:basetriple}. Now let $k > 2$ and $|E(H)| > 0$; let $e = uv \in E(H)$ and suppose that $H \setminus \sset{e}$ 
is $(m_1, k)$-forcible and that $H$ is $(m_2, k-1)$-forcible. Let $m = m_1^{m_2}$. We claim that $H$ is $(m, k)$-forcible.

For an integer $s\ge 1$, we say an
\emph{$s$-star} is a graph $J$ with $s(|H|-1) + 1$ vertices partitioned 
into sets $V_1, \dots, V_s, \sset{w}$ such that
\begin{itemize}
\item $|V_i| = |H|-1$ for all $i \in \sset{1, \dots, s}$;
\item $V_i$ is anticomplete to $V_j$ for all $i, j \in \sset{1, \dots, s}$ with $i \neq j$; and
\item $J[V_i \cup \sset{w}]$ is isomorphic to $H\setminus \sset{e}$ and $w$ maps to $u$ under this isomorphism, for all $i \in \sset{1, \dots, s}$. 
\end{itemize}
The vertex $w$ is called the \emph{centre} of $J$. 
\\
\\
(1) {\em For every $s \geq 1$, the $s$-star is $(m_1^s, k)$-forcible.}
\\
\\
We prove this by induction on $s$. For $s = 1$, this follows since $H \setminus \sset{e}$ is $(m_1, k)$-forcible and 
is isomorphic to the $1$-star. Now let $s > 1$. By induction, the $(s-1)$-star is $(m_1^{s-1}, k)$-forcible, and so 
by \ref{lem:union}, the disjoint union of $m_1$ graphs, each an $(s-1)$-star, 
is $(m_1^s, k)$-forcible. Starting with a graph $G_0$ with $m_1^s$ vertices and no edges, player A uses the strategy that forces
this disjoint union, until at the end of round $p$, say, $G_{p}$ has an induced subgraph with $m_1$ components $H_1\ll H_m$,
each an $(s-1)$-star.
For $i \in \sset{1, \dots, m_1}$, let $u_i$ denote the centre of $H_i$, and let
$U = \sset{u_1\ll u_{m_1}}$. It follows that $|U| = m_1$ and $U$ is stable in $G_p$. By applying the strategy 
for $H \setminus \sset{e}$ to $U$ starting at round $p+1$, it follows that there is a strategy for player A that 
produces at the end of some round $q$ a 
graph $G_q$ containing an induced subgraph that consists of 
$m_1$ disjoint $(s-1)$-stars, pairwise anticomplete except for edges with both ends in a set $U$ as defined above, 
and a subset $U'\subseteq U$ such that $G_q[U']$ is isomorphic to $H \setminus \sset{e}$.
Let $i \in \sset{1, \dots, m_1}$ such that $u_i \in U'$ 
maps to $u$ under the isomorphism 
between $G_p[U']$ and the graph $H \setminus \sset{e}$. It follows that $U' \setminus \sset{u_i}$ is anticomplete 
to $V(H_i) \setminus \sset{u_i}$, and therefore $G_p[U' \cup V(H_i)]$ is an $s$-star. This proves (1). 

  \bigskip

By (1), it follows that an $m_2$-star is $(m_1^{m_2}, k)$-forcible. This implies that in the $k$-tuple game on $m$ vertices, 
player A can guarantee that in some round $s$, $G_s$ contains an $m_2$-star $H'$ with vertex set 
$V_1 \cup \dots \cup V_{m_2} \cup \sset{w}$ (with notation as before). Let $w$ be the centre of $H'$, and for 
$i \in \sset{1, \dots, m_2}$, let $v_i$ be the vertex corresponding to $v$ in the isomorphism between 
$G_s[V_i \cup \sset{w}]$ and $H \setminus \sset{e}$. Let $V = \sset{v_1\ll v_{m_2}}$. By the definition of an $m_2$-star, it follows 
that $V$ is a stable set in $G_s$, and $|V| = m_2$. Now player A uses the winning strategy of the $(k-1)$-tuple game for $H$
on $m_2$ vertices by starting with $G_s[V]$; except in every round, player A picks $T \cup \sset{w}$ instead of the 
set $T \subseteq V$ of $(k-1)$ vertices that the strategy for the $(k-1)$-tuple game produces. If in round $s' > s$, player B 
adds an edge incident with $w$, say $wv_i$, then $G_{s'}[V_i \cup \sset{w}]$ is isomorphic to $H$, and the result follows. 
Therefore, we may assume that in every round $s' > s$, player A picks a set $T \cup \sset{w}$, and player B adds at least 
one edge with both ends in $T$. Since $H$ is $(m_2, k-1)$-forcible, it follows that for some $s' > s$, $G_{s'}[V]$ 
contains $H$. This proves \ref{lem:triplemain}.~\bbox

Let us digress for a moment. For a graph $H$ and an integer $k\ge 2$, let $m(H,k)$ be the smallest $m$ such that
$H$ is $(m,k)$-forcible. If $H$ is a complete graph $K_t$ then
$m(H,k)$ equals the Ramsey number $r(t,k)$, the smallest integer $r$ such that every graph with at least $r$ vertices has either
a clique of size $t$ or a stable set of size $k$. (To see this, if we are playing on $m\ge r(t,k)$ vertices, then player A must win, with arbitrary play; 
and if $m<r(t,k)$, then player B has a winning strategy, since
there is a graph $G$ on this vertex set with no clique of size $k$ and no stable set of size $t$, and in each turn player B fills in
the edges of $G$ that lie within the stable set passed by player A.)
Thus in this case $m(H,k)$ is single-exponential in $k$. In general,
\ref{lem:triplemain} shows that $m(H,k)$ exists, and gives an upper bound on $m(H,k)$
that is ``tower-type'' in $k$. This can be reduced to something doubly-exponential in $k$, by improving (1) above, 
showing that every $s$-star is $(M, k)$-forcible where $M$ is some function that is polynomial in $m_1$ and $s$.
We could show this as follows (sketch). 

Take an $m_1$-uniform hypergraph $J$ such that
\begin{itemize}
\item  $J$ has minimum degree at least $m_1s$;
\item the bipartite graph with bipartition $(E(J), V(J))$ associated with $J$ has no cycles of length less than eight; and
\item $|V(J)|$ is at most polynomial in $m_1,s$.
\end{itemize}
(One can show that such a hypergraph exists.) Now play the game within each hyperedge of $J$, and one easily
obtains an $s$-star.

We have not determined in general whether $m(H,k)$ is singly- or doubly-exponential in $k$.

\section{Sparse $k$-tuples} \label{sec:sparsek}

In this section we prove a lemma that is used for all the difficult results of the paper, \ref{lem:splusone} below. Its proof uses the
game from the previous section, and some other preliminaries.

\begin{thm} \label{lem:subsat}
Let $G, H$ be graphs, and let $\alpha > 0$ such that $G$ is $(\alpha, H)$-saturated.
Let $H'$ be an induced subgraph of $H$. Then $G$ is $(\alpha, H')$-saturated.
\end{thm}
\Proof
Let $n = |G|$, and let $\mathcal{T}$ be the set of copies of $H'$ in $G$. For every copy $\phi$ of $H$ in $G$, 
$\phi|_{V(H')}$ is a copy of $H'$ in $G$; we say that $\phi$ \emph{came from}  
$\phi|_{V(H')}$. For every $\phi' \in \mathcal{T}$, there are at most $n^{|H| - |H'|}$ copies of $H$ 
that came from $\phi'$ (since there are at most $n^{|H| - |H'|}$ ways to extend $\phi'$ from a function 
from $V(H')$ to $V(G)$, to a function from $V(H)$ to $V(G)$). It follows that 
$$|\mathcal{T}| \geq n^{{|H'| - |H|}} \alpha n^{|H|} = \alpha n^{|H'|},$$ 
and so $G$ is $(\alpha, H')$-saturated. This proves \ref{lem:subsat}.~\bbox

Let $G$ be a graph with $n$ vertices, and let $A_1, \dots, A_k \subseteq V(G)$ be pairwise disjoint such that
\begin{itemize}
\item for all $i \in \sset{1, \dots, k}$, $|A_i| \geq \alpha n$; and
\item for all $i, j \in \sset{1, \dots, k}$ with $i \neq j$, $(A_i, A_j)$ is a $c$-sparse pair. 
\end{itemize}
Then we call $A_1, \dots, A_k$ a {\em $c$-sparse $(\alpha, k)$-tuple} in $G$.
Thus a $c$-sparse $(x,x)$-pair is a $c$-sparse $(x/n,2)$-tuple. 

\begin{thm} \label{lem:copies}
Let $G, H$ be graphs. Let $S \subseteq V(H)$ be a stable set with $|S| = k\ge 2$; and let $0\le \alpha, c \le 1$ 
such that $G$ is $\left(\alpha, H\right)$-saturated. Then either $G$ contains a $c$-sparse 
$\left(\frac{\alpha}{2 k^k}, k\right)$-tuple, or there is an $S$-successor $H'$ of $H$ such that $G$ is 
$\left(\frac{c\alpha }{2^{k^2} k^k}, H'\right)$-saturated.
\end{thm}
\Proof
Let $\ell = |H|$ and $n=|G|$.
For a copy $\phi$ of $H$, we say that $\phi|_{H \setminus S}$ is 
the \emph{anchor} of $\phi$. 
We say that a copy $\psi$ of $H\setminus S$ is \emph{weighty} if there are at least $\frac{1}{2} \alpha n^k$ 
copies of $H$ in $G$ with anchor $\psi$. Let $\mathcal{T}$ be the set of weighty copies of $H\setminus S$. 
Since there are at most $n^{\ell-k}$ copies of $H\setminus S$, there are 
at most $\frac{1}{2} \alpha n^\ell$ copies of $H$ whose anchors are not weighty. Consequently, there are at least 
$\frac{1}{2} \alpha n^\ell$ copies of $H$ whose anchors are weighty.
For $\psi \in \mathcal{T}$, let $\mathcal{A}(\psi)$ be the set of copies of $H$ with anchor $\psi$. 
It follows that 
$$\sum_{\psi \in \mathcal{T}} |\mathcal{A}(\psi)| \geq  \alpha n^\ell/2.$$

Let
$S = \sset{s_1, \dots, s_k}$.
For $i \in \sset{1, \dots, k}$, we let $U_i = \sset{\phi(s_i) : \phi \in \mathcal{A}(\psi)}$. Now let $V_1, \dots, V_k$ 
be a random partition of $U = U_1 \cup \dots \cup U_k$ in which, for all $i \in \sset{1, \dots, k}$,
every vertex of $U$ is in $V_i$ with probability $1/k$ 
independently.  Let $\phi \in  \mathcal{A}(\psi)$. It follows that the probability 
that $\phi(s_i) \in V_i$ for all $i \in \sset{1, \dots, k}$ is $1/k^k$. Therefore there is a choice of $V_1, \dots, V_k$ 
such that 
$$|\sset{\phi \in  \mathcal{A}(\psi) : \phi(s_i) \in V_i \text{ for all } i \in \sset{1, \dots, k}}| \geq |\mathcal{A}(\psi)|/k^k.$$
Fix such a choice of $V_1, \dots, V_k$,
and let $W_i = U_i \cap V_i$ for all $i \in \sset{1, \dots, k}$. 
It follows that  
$$|W_1|\cdots |W_k|\ge \sset{\phi \in  \mathcal{A}(\psi) : \phi(s_i) \in W_i \text{ for all } i \in \sset{1, \dots, k}}| \geq |\mathcal{A}(\psi)|/k^k,$$
since $\phi(s_i) \in U_i$ for all $\phi \in \mathcal{A}(\psi)$ and $i \in \sset{1, \dots, k}$.
Since  
$|\mathcal{A}(\psi)|/k^k \ge \alpha n^kk^{-k}/2,$ and $|W_1|\cdots |W_k|\le n^{k-1}|W_i|$,
it follows that 
$|W_i| \geq \frac{\alpha}{2k^k} n$ for all $i \in \sset{1, \dots, k}$. If for all distinct $i, j \in \sset{1, \dots, k}$,
$(W_i, W_j)$ is a $c$-sparse pair, then $W_1, \dots, W_k$ is a $c$-sparse 
$\left(\frac{c\alpha}{2 k^k}, k\right)$-tuple, and \ref{lem:copies} follows. Therefore, we may assume that there exist 
distinct $i, j \in \sset{1, \dots, k}$ such that $(W_i, W_j)$ is not $c$-sparse. 

Let $\mathcal{W}(\psi)$ be the set of all $w = (w_1, \dots, w_k)$ such that $w_{h} \in W_{h}$ for all $h \in \sset{1, \dots, k}$, 
and $w_i, w_j$ are adjacent.
For $w = (w_1, \dots, w_k)\in \mathcal{W}(\psi)$, let
$\phi_w(v) = \psi(v)$ if $v \in V(H) \setminus S$, and $\phi_w(s_{h}) = w_{h}$ for all $h \in \sset{1, \dots, k}$.
Since $w_h\in U_h$ for each $i$, and $w_i, w_j$ are adjacent, it follows that
$\phi_w$ is a copy of an $S$-successor of $H$.

Since the pair $(W_i, W_j)$ is not $c$-sparse, 
$$|\mathcal{W}(\psi)|\ge c|W_1| \cdots |W_k| \ge c|\mathcal{A}(\psi)|/k^k.$$
But $\sum_{\psi \in \mathcal{T}} |\mathcal{A}(\psi)| \geq \frac{1}{2} \alpha n^\ell,$ 
and so
$$\sum_{\psi \in \mathcal{T}} |\mathcal{W}(\psi)| \geq \frac{c\alpha}{2 k^k} n^\ell.$$ 
This implies that 
$G$ contains at least $\frac{c\alpha}{2 k^k} n^\ell$ copies of $S$-successors of $H$. Now $H$
has at most $2^{k^2-1}$ distinct $S$-successors, since that bounds the number of distinct graphs on $k$ vertices;
and therefore, there is an $S$-successor $H'$ of $H$ such 
that $G$ contains at least $2^{-k^2} k^{-k} \alpha cn^\ell$ copies of $H'$. It follows that $G$ 
is $\left(\frac{c\alpha}{2^{k^2} k^{k}}, H'\right)$-saturated. This proves \ref{lem:copies}.~\bbox

\begin{thm} \label{lem:stablesat}
Let $H$ be a graph with $|E(H)| = 0$, and let $\vare >0$ with $\vare|H| \leq 1/2$. Let $G$ be an $\vare$-bounded graph.
Then $G$ is $\left(2^{-|H|}, H\right)$-saturated.
\end{thm}
\Proof
We prove this by induction on $|H|$. For $|H| = 0$, $H$ is the null graph, and so $G$ is trivially $(1, H)$-saturated. 
Now let $|H| > 0$, and let $v \in V(H)$. From the inductive hypothesis, it follows that $G$ is 
$\left(2^{-|H|+1}, H \setminus \sset{v}\right)$-saturated. Let $\phi$ be a copy of $H \setminus \sset{v}$ in $G$. 
Since $G$ is $\vare$-bounded, it follows that there are at most $\vare (|H| -1)|G|\le |G|/2$ vertices of $G$ that are equal to or adjacent to
a vertex in the image of $\phi$; and so there are at least $|G|/2$ that are not. Consequently there are at least $|G|/2$
vertices $w$ such that the function $\phi_w$ is a copy of $H$, where $\phi_w(x) = \phi(x)$ if $x \in V(H) \setminus \sset{v}$, 
and $\phi(v) = w$. Summing over $\phi$, it follows that 
there are at least $2^{-|H| } |G|^{|H|}$ 
copies of $H$ in $G$, and so $G$ is $\left(2^{-|H|}, H\right)$-saturated. This proves \ref{lem:stablesat}.~\bbox

\begin{thm} \label{lem:splusone}
Let $k\ge 2$ and $t\ge 0$ be integers and let $H$ be a graph.
Then there exist $S, \vare > 0$ such that for
every $\vare$-bounded graph $G$, and for all $c$ with $0< c\le 1$, either $G$ is $\left(\vare c^S, H\right)$-saturated, or $G$
contains a $c^{s+t}$-sparse $(\vare c^s, k)$-tuple for some $s \in \sset{0, \dots, S}$.
\end{thm}

\Proof Let $\ell=|H|$; we may assume that $\ell>0$.
By \ref{lem:triplemain}, it follows that there exists $L \in \mathbb{N}$ such that $H$ is $(L, k)$-forcible, and so $L\ge \ell$.
Let $r = {L \choose 2}+1$. For $i \in \sset{0, \dots, r}$, let 
$\delta_i = \left(\frac{1}{k^k 2^{k^2}}\right)^i 2^{-L}$; and let $\vare = \min(1/(2L), \delta_r)$.
Then $\delta_{i+1} \leq \frac{1}{k^k 2^{k^2}} \delta_i$ for each $i$.
Let $S=2^{r+t}$ and $s_i= 2^{i+t}-t$ for each $i$.
Let $G$ be $\vare$-bounded, and let $0< c\le 1$. We will show that either $G$ is $\left(\vare c^S, H\right)$-saturated, or $G$
contains a $c^{s+t}$-sparse $(\vare c^s, k)$-tuple for some $s \in \sset{0, \dots, S}$. Thus we may assume:
\\
\\
(1) {\em For $0\le i< r$, $G$ does not contains a $c^{s_i+t}$-sparse $(\vare c^{s_i}, k)$-tuple.}
\\
\\
Let $H_0$ be an $L$-vertex graph with no edges. By  \ref{lem:stablesat}, it follows that $G$ is $\left(\delta_0, H_0\right)$-saturated. 
We will play the $k$-tuple game for $H$ on $L$ vertices starting with the graph $H_0$; and the graph passed to player A at the 
start of round $i$ is denoted by $H_i$. Since every round (except the last) adds an edge, it follows that there are at most $r$ rounds. Player A
will use an optimal strategy, one that will guarantee that at some round, $H_i$ will contain $H$. 
Player B will use the following strategy. In round $i$, when player B is presented with a stable subset $T$ of $H_{i-1}$ of cardinality $k$, player B
will return a graph $H_i$ determined as follows. By 
\ref{lem:copies}  applied with 
$\alpha=\delta c^r$, $m = t+s_{i-1}$,  $\delta = \delta_{i-1}$, $c=c^m$, and $s =s_{i-1}$, it follows that either $G$ has a
$c^{t+s_{i-1}}$-sparse $(\delta_i c^{s_{i-1}}, k)$-tuple, contrary to (1), or, 
since $s_i = 2s_{i-1}+ t$, there is
a $T$-successor $H_i$ of $H_{i-1}$ such that $G$ is $\left(\delta_i c^{s_i}, H_i\right)$-saturated. Then player B returns $H_i$.

It follows that at the start of each round $i$ of the game, $G$ is $\left(\delta_{i-1} c^{s_{i-1}}, H_{i-1}\right)$-saturated.
If player A wins the game in this round, then $H_{i-1}$ contains $H$, and $G$ is 
$\left(\delta_{i-1} c^{s_{i-1}}, H_{i-1}\right)$-saturated. Since $H$ is an induced subgraph of $H_{i-1}$, it follows that $G$ is 
$\left(\vare c^S, H\right)$-saturated by  \ref{lem:subsat}. Since player A wins in some round,
this proves \ref{lem:splusone}.~\bbox

We remark that in the case of $k=2$, we can bound $s$ by $|E(H)|$:

\begin{thm}\label{2csparse}
Let $k \in \mathbb{N}$ and let $H$ be a graph with $s$ edges. Then there exists $\vare > 0$ such that 
for every $\vare$-bounded graph $G$, and for all $c$ with $0\le c\le 1$, either $G$ is  $\left(\vare c^s, H\right)$-saturated, or $G$ contains 
a $c$-sparse $(\vare c^s, 2)$-tuple.
\end{thm}
\Proof
This follows by changing slightly the strategy of player B in the proof of \ref{lem:splusone}: just maintain that 
$G$ is $\left(\delta_0 c^i, H_i\right)$-saturated. The result follows since $H$ is $(|H|, 2)$-forcible and player A has a strategy 
that forces $H$ in $s$ rounds.~\bbox

A special case of \ref{lem:splusone} is of interest and worth stating separately:
\begin{thm}\label{lem:kcsparse}
Let $k \ge 2$ be an integer, and let $H$ be a graph. Then there exist $s \in \mathbb{N}$ and $\vare> 0$ such that for 
every $\vare$-bounded graph $G$, and for all $c$ with $0\le c\le 1$, either $G$ is  $\left(\vare c^s, H\right)$-saturated, or $G$ 
contains a $c$-sparse $(\vare c^s, k)$-tuple.
\end{thm}
\Proof
By \ref{lem:splusone} with $t=0$, there exist $S, \vare > 0$ such that for
every $\vare$-bounded graph $G$, and for all $c$ with $0\le c\le 1$, either $G$ is $\left(\vare c^S, H\right)$-saturated,
or $G$
contains a $c^{r}$-sparse $(\vare c^r, k)$-tuple for some $r \in \sset{0, \dots, S}$. In the second case
$G$
contains a $c$-sparse $(\vare c^S, k)$-tuple, so in both cases the theorem holds taking $s=S$.~\bbox

\section{Excluding general graphs}\label{sec:ehpimprove}

The results of the previous section can be applied to deduce several results about excluding general graphs, 
that we obtain in this section.
We need first:
\begin{thm}\label{getsparse}
Let $G$ be a graph with $n$ vertices and at most $\vare n(n-1)/2$ edges with $n\ge 2\vare^{-1}$.
Then there is an induced subgraph $J$ with $|J|\ge n/2$
such that $J$ is $2\vare$-bounded.
\end{thm}
\Proof
Choose distinct $v_1\ll v_k \in V(G)$ with $k$ maximum such that for $1\le i\le k$,
$v_i$ has at least $2\vare (n-i+1)-1$ neighbours in $V(G)\setminus \{v_1\ll v_i\}$.
Let $m=\ceil{n/2}$.
If $k\ge m$, then there are at least
$$\sum_{1\le i\le m} (2\vare (n-i+1)-1)=2\vare m(n-(m-1)/2)-m\ge 3\vare mn/2-m\ge (3\vare n/2-1)n/2$$
edges in $G$ with an end in $\{v_1\ll v_m\}$.
Consequently $(3\vare n/2-1)n/2\le \vare n(n-1)/2$, contradicting that
$n\ge 2\epsilon^{-1}$. So $k\le n/2$. But from the maximality of $k$,
$G[V(G)\setminus \{v_1\ll v_k\}]$ is $2\vare$-bounded, and therefore satisfies the theorem.~\bbox

We use \ref{getsparse} to prove a consequence of \ref{thm:copies}.
\begin{thm} \label{cor:rodlcount}
For every graph $H$ and every $\vare > 0$ there exist $\alpha, \delta > 0$ such that for every graph $G$, if $G$ is
not $(\alpha, H)$-saturated, then either $|G| \leq 4(\delta\vare)^{-1}$ or $G$ contains an induced subgraph $J$
with $|J| \geq \delta |G|$ such that one of $J,\overline{J}$ is $\vare$-bounded.
\end{thm}
\Proof
Let $H$ be a graph, and let $\vare > 0$. Let $\alpha, \delta > 0$ satisfy \ref{thm:copies}, with $\vare,\delta$ replaced by
$\vare/2,2\delta$ respectively.
Now let $G$ be a graph that is not $(\alpha, H)$-saturated. By \ref{thm:copies},
$G$ contains an induced subgraph $J$ with $|J| \geq 2\delta |G|$ and such that
either $|E(J)| \leq \vare |J|(|J| -1)/4$ or $|E(\overline{J})| \leq \vare |J|(|J| -1)/4$. In the first case, by \ref{getsparse},
either $|J|< 4\vare^{-1}$ and hence $|G|< 4(\vare\delta)^{-1}$, or $G$ contains an $\vare$-bounded induced subgraph with at
least $|J|/2 \geq \delta |G|$ vertices. In the second case we use the same argument in the complement. This proves \ref{cor:rodlcount}.~\bbox

This is used to prove \ref{EHdensityimprove}. Alternatively, \ref{EHdensityimprove} can be deduced 
from corollary 3.3 of~\cite{foxsudakov} (as written,
this shows that $G$ is not $H$-free, but the same
proof can be used to show that $G$ is $(\alpha, H)$-saturated). We restate \ref{EHdensityimprove}:

\begin{thm}\label{restate}
For every graph $H$ and $k \in \mathbb{N}$, there exist $\vare, s, K >0$ such that for every graph $G$ with $n>K$ vertices,
and every $c$ with $0\le c\le 1$, if $G$ is not $(\vare c^s, H)$-saturated, then there are pairwise disjoint subsets $A_1, \dots, A_k\subseteq V(G)$ such that
either:
\begin{itemize}
\item $(A_i, A_j)$
is a $c$-sparse $(\vare c^s n, \vare c^s n)$-pair for all distinct $i, j \in \sset{1, \dots, k}$; or
\item $(A_i, A_j)$
is a $(1-c)$-dense $(\vare c^s n, \vare c^s n)$-pair for all distinct $i, j \in \sset{1, \dots, k}$.
\end{itemize}
\end{thm}
\Proof
Let $\vare', s$ satisfy \ref{lem:kcsparse} both for $H$ and
for $\overline{H}$. Let $\alpha, \delta\le 1$ be as in \ref{cor:rodlcount} for $H$ and $\vare'$. Let $K = 4(\delta\vare')^{-1}$.
Let
$\vare = \min\left(\alpha, \vare' \delta^{|H|}\right).$
We claim that $s, K, \vare$ satisfy the theorem.

Let $G$ be a graph. By \ref{cor:rodlcount}, it follows that either $G$ is $(\alpha, H)$-saturated (and thus $(\vare, H)$-saturated),
or $|G| \leq 4(\delta\vare')^{-1} = K$, or $G$ contains an induced subgraph $J$ with $|J| \geq \delta |G|$ such that either $J$ or
$\overline{J}$ is $\vare'$-bounded. We may assume the third of these holds. Let $0\le c\le 1$.

Suppose first that $J$ is $\vare'$-bounded. By \ref{lem:kcsparse}, it follows that either $J$ is
$\left(\vare' c^{s}, H\right)$-saturated
(and so $G$ is $\left(\vare' \delta^{|H|} c^s, H\right)$-saturated), or $J$ contains a $c$-sparse $(\vare' c^{s}, k)$-tuple.
We may
assume the latter; but then $G$ contains a $c$-sparse $(\vare' \delta^{|H|} c^s, k)$-tuple, and \ref{restate} follows.
In the case when $\overline{J}$ is $\vare'$-bounded, we apply the same argument in the complement, using $\overline{H}$
instead of $H$. This proves \ref{restate}.~\bbox

Before the next result we need two easy lemmas:

\begin{thm}\label{getanticomp1}
Let $k\in \mathbb{N}$, let $G$ be a graph, and let $P_1\ll P_k$ be pairwise disjoint subsets of $V(G)$, 
where $|P_i|=p_i$ for $1\le i\le k$. Let $d_1\ll d_k$ be such that for all distinct $i,j\in \{1\ll k\}$,
and every $v\in P_i$, $v$ has at most $d_j$ neighbours in $P_j$. Let $q\in \mathbb{N}$ with $((k-1)d_i+1)q\le p_i$
for $1\le i\le k$. Then there
are subsets $Q_i\subseteq P_i$ for $1\le i\le k$, each of cardinality $q$, and pairwise anticomplete.
\end{thm}
\Proof
For $k=0$ the result is vacuously true, so we assume that $k\ge 1$ and that the result holds for $k-1$. Choose
$Q_k\subseteq P_k$ of cardinality $q$ (this is possible since $q\le ((k-1)d_i+1)q\le p_i$), and for $1\le i\le k-1$ 
let $P_i'$ be the set of vertices in $P_i$ with no neighbour in $Q_k$.
Thus $|P_i'|\ge p_i-qd_i\ge ((k-2)d_i+1)q$ for $1\le i\le k-1$, so from the inductive hypothesis 
there exist $Q_i\subseteq P_i'$ of cardinality $q$ 
for $1\le i\le k-1$,
pairwise anticomplete; and they are all anticomplete to $Q_k$. This proves \ref{getanticomp1}.~\bbox

This extends to:
\begin{thm}\label{getanticomp2}
Let $k\in \mathbb{N}$, let $G$ be a graph, let $0\le c\le 1$, and let $P_1\ll P_k$ be disjoint subsets of $V(G)$, pairwise $c$-sparse.
Let $|P_i|=p_i$ for $1\le i\le k$.
Let $q\in \mathbb{N}$ such that $2q(2(k-1)^2cp_i+1)\le p_i$ for $1\le i\le k$.
Then there exist $Q_i\subseteq P_i$ for $1\le i\le k$, each of cardinality $q$, and pairwise anticomplete.
\end{thm}
\Proof We may assume that $c>0$.
For all distinct $i,j\in \{1\ll k\}$, let $B_{i,j}$ be the set of vertices $v\in P_i$ with more than $2(k-1)cp_j$ neighbours in $B_j$.
Since there are at most $cp_ip_j$ edges between $P_i$ and $P_j$, it follows that $|B_{i,j}|\le cp_ip_j/(2(k-1)cp_j)=p_i/(2(k-1))$.
By taking the union of $B_{i,j}$ for all $j\ne i$, we deduce that there are at most $p_i/2$ vertices in $P_i$ that have
more than $2(k-1)cp_j$ neighbours in $B_j$ for some $j\ne i$; and so there are at least $p_i/2$ vertices in $P_i$ that have at most
$2(k-1)cp_j$ neighbours in $P_j$ for each $j\ne i$. By \ref{getanticomp1}, there 
are subsets $Q_i\subseteq P_i$ for $1\le i\le k$, each of cardinality $q$, and pairwise anticomplete. This proves \ref{getanticomp2}.~\bbox

The following result shows that there are $k$ sets, each of size $\vare n^\vare$, and pairwise anticomplete, if we exclude a graph $H$ 
as an induced subgraph of an $\vare$-bounded graph. This is similar to \ref{EHP}, except that we assume sparsity, and guarantee 
anticomplete sets, and we get more than two anticomplete sets. 

\begin{thm} \label{lem:kpoly2}
Let $H$ be a graph and $k\ge 2$ be an integer. Then there exists $\vare > 0$ such that in every  $\vare$-bounded $H$-free 
graph $G$ with $|G| = n\ge 2$, there are $k$ disjoint subsets of $V(G)$, pairwise anticomplete and each of cardinality at 
least $\vare n^{\vare}$.
\end{thm}
\Proof
Let $\vare', s > 0$ be as in  \ref{lem:kcsparse}.
Let
$$\vare = \min\left(\vare'/8, \frac{1}{2|H|}, \frac{1}{16k^2}, \frac{1}{s+1}\right),$$ 
and let $G$ 
be $\vare$-bounded and $H$-free, and let $n=|G|$. By  \ref{lem:stablesat}, it follows that $G$ contains a stable set of size $k$; therefore, we may 
assume that $\vare n^{\vare} > 1$. Let $c = n^{-1/(s+1)}$.

By  \ref{lem:kcsparse}, it follows that either $G$ is $\left(\vare' c^s, H\right)$-saturated or $G$ contains a
$c$-sparse $(\vare' c^s, k)$-tuple $A_1, \dots, A_k$, and since $G$ is $H$-free, the latter holds. 
Let $q=\ceil{\vare n^{\vare}}$. By \ref{getanticomp2} it suffices to show that
$2q(2(k-1)^2cp_i+1)\le p_i$ for $1\le i\le k$, where $p_i=|A_i|$. Thus, it suffices to check that
$2q(2(k-1)^2cp_i)\le p_i/2$ and $2q\le p_i/2$, that is, 
$8qc(k-1)^2\le 1$ and $4q\le \vare' c^sn$. Since $\vare n^{\vare} > 1$, it follows that $q\le 2\vare n^{\vare}$;
so it suffices to show that  $16\vare n^{\vare}c(k-1)^2\le 1$ and $8\vare n^{\vare}\le \vare' c^sn$.

For the first, since $c = n^{-1/(s+1)}$, we must show that
$16\vare n^{\vare} n^{-1/(s+1)}(k-1)^2\le 1$, and this is true since $\vare\le 1/(s+1)$ and $16\vare(k-1)^2\le 1$.
For the second, we must show that
$8\vare n^{\vare}\le \vare' n^{1-s/(s+1)}$, and this is true since $\vare\le 1/(s+1)$ and $\vare\le \vare'/8$.
This proves \ref{lem:kpoly2}.~\bbox

The next result is an improvement of \ref{EHP} in the $\vare$-bounded case. 
\begin{thm} \label{lem:onepluseps}
Let $H$ be a graph. Then there exists $\vare > 0$ such that if $G$ is $H$-free and $\vare$-bounded, then $G$ has an
anticomplete pair $(A, B)$ with $|A|\cdot |B| \geq \vare n^{1+\vare}$.
\end{thm}

\Proof
Let $S, \vare'$ be as in  \ref{lem:splusone}, setting $k = 2$ and $t=1$. 
We may assume that $\vare'\le 1/4$.
Let 
$$\vare=\min\left((\vare')^2/2, \frac{1}{2S+1}\right);$$ 
we claim that $\vare$ satisfies the theorem. Let $G$ be $\vare$-bounded, and let $|G| = n$. 
Let $c = n^{-1/(2S+1)}$.
It follows from  \ref{lem:splusone} that either $G$ is
$\left(\vare' c^S, H\right)$-saturated, or $G$ contains a $c^{s+1}$-sparse $(\vare' c^s, 2)$-tuple for some
$s \in \sset{0, \dots, S}$. Since $G$ is $H$-free, the latter holds, and so $G$ contains a $c^{s+1}$-sparse
$(\vare' c^s n, \vare' c^sn)$-pair $(A, B)$ for some $s \in \sset{0, \dots, S}$. Let $t=(s+1)/(2S+1)$, and $m=\ceil{\vare' n^{t}}$. 
\\
\\
(1) {\em We may assume that $\vare' n^{t}\ge 1/2$, and $\vare' n^{t}\le |A|$, and $\vare' n^{t}c^{s+1}\le 1/4$.}
\\
\\
Let $v\in V(G)$; then we may assume that $|\sset{v}|\cdot |V(G) \setminus N[v]| <\vare n^{(1+\vare)}$, for otherwise the theorem holds.
But $v$ has at least $(1-\vare)n$ non-neighbours, so $(1-\vare)n<\vare n^{1+\vare}$, and hence 
$$n^{\vare}>1/\vare-1\ge 1/\vare'-1\ge 1/(2\vare').$$
Consequently $\vare' n^{t}\ge \vare' n^{\vare}\ge 1/2$,
so the first holds.
The second holds since $n^{t}\le c^s n$,
and the third since $n^{t}c^{s+1}=1$. This proves (1).

\bigskip

By \ref{sparsetoanti}, taking $x=\vare' n^{t}$ and with $c$ replaced by $c^{s+1}$, 
there is an anticomplete $(\vare' n^{t}, |B|/2)$ pair $(A', B')$ with $A'\subseteq A$ and $B'\subseteq B$.
Then 
$$|A'|\cdot|B'|\ge \left(\vare' n^{t}\right)\left(\frac12 \vare'  c^sn\right)= \frac12 \vare'^2 n^{1+1/(2S+1)}\ge \vare n^{1+\vare}.$$
This proves \ref{lem:onepluseps}.~\bbox

This implies \ref{easyEH}, which we restate.

\begin{thm} \label{lem:onepluseps2}
For every graph $H$ there exists $\vare >0$ such that in every $H$-free graph $G$ with $n>1$ vertices,
there is a pure pair $(A,B)$
with $|A|\cdot|B|\ge \vare n^{1+\vare}$.
\end{thm}
\Proof
Let $\vare'$ satisfy \ref{lem:onepluseps} for both $H$ and $\overline{H}$.
Let $\delta$ satisfy \ref{cor:rodlcount} for $\vare'$ and $H$. Let
$$\vare=\min\left(\vare'\delta^{1+\vare'},  (\delta\vare')^2/16\right).$$
Let $G$ be $H$-free, and let $n = |G|\ge 2$.  By \ref{cor:rodlcount}, it follows that either $n \leq 4/(\delta\vare')$, or
$G$ contains an
induced subgraph $J$ with at least $\delta n$ vertices such that one of $J$, $\overline{J}$ is $\vare'$-bounded.
If $n \leq 4/(\delta\vare')$, then $n^{1+\vare}\le n^2\le \vare^{-1}$ from the definition of $\vare$.
Choose distinct $u,v\in V(G)$, and then $\sset{u}, \sset{v}$ is an
pure pair $(A, B)$ with $|A|\cdot |B| =1\geq \vare n^{1+\vare}$, and \ref{lem:onepluseps2} holds.
Therefore, we may assume that $G$ contains an induced subgraph $J$ with at least $\delta n$ vertices
such that one of $J$, $\overline{J}$ is $\vare'$-bounded.
If $J$ is $\vare'$-bounded, then \ref{lem:onepluseps} implies that $J$
contains an anticomplete pair $(A, B)$ with
$$|A|\cdot |B| \geq \vare' |J|^{1+\vare'} \geq \vare n^{1+\vare}.$$
If $\overline{J}$ is $\vare'$-bounded, we apply the same argument in the complement, using $\overline{H}$, obtaining
a complete pair in $G$. This proves \ref{lem:onepluseps2}.~\bbox

We have given a long and complicated proof for \ref{lem:onepluseps2}, since it is a consequence of other results
that we needed anyway; but \ref{lem:onepluseps2} can be proved directly, much more easily, using a minor variant of the original proof of \ref{EHP}
by Erd\H{o}s, Hajnal and Pach~\cite{EHP}, as follows. We use
the following lemma.
For $k\ge 1$, define $e_k = 1-2^{1-k}$.
\begin{thm}\label{easylemma}
Let $H$ be a graph with $k\ge 1$ vertices $h_1\ll h_k$, and let $t\ge 5^{2^{k-2}}$ be a real number.
Let $G$ be a $k$-partite graph, with parts $V_1\ll V_k$, each of cardinality at least $5t^{e_k}$.
Then either
\begin{itemize}
\item for $1\le i\le k$ there exists $v_i\in V_i$ such that for $1\le i<j\le k$, $v_i,v_j$ are adjacent in $G$
if and only if $h_i, h_j$ are adjacent in $H$; or
\item there exist $i,j$ with $1\le i<j\le k$ and subsets $A\subseteq V_i$ and $B\subseteq V_j$  such that
$A,B$ are complete or anticomplete to each other, and $|A|\cdot|B|\ge t$.
\end{itemize}
\end{thm}
\Proof We may assume that $H$ is a complete graph, by replacing all edges between $V_i, V_j$ by the bipartite complement
if $h_i, h_j$ are nonadjacent.
If $k=1$ the result is trivial. We assume $k>1$ and proceed by induction on $k$.

Define $n= 5t^{e_k} $ and $d = 5t^{e_{k-1}}$.
If there exists $v_1\in V_1$ such that $v_1$ has at least $d$ neighbours in
each of $V_2\ll V_k$, then
the result follows by induction (applied to $H\setminus \{h_1\}$ and the sets $N[v_1]\cap V_i\;(2\le i\le k)$), since $t^{2^{2-k}}\ge t^{2^{1-k}}\ge 5$.

So we may assume that each vertex in $V_1$ has fewer than $d$ neighbours in one of $V_2\ll V_k$; and so we may assume that
at least $n/(k-1)$ vertices in $V_1$ have fewer than $d$ neighbours in $V_2$. Now
since $t\ge 5^{2^{k-2}}$ by hypothesis, it follows that
$$t^{e_{k-1}}\ge 5^{2^{k-2}(1-2^{2-k})}= 5^{2^{k-2}-1}\ge k-1.$$ Consequently
$n/(k-1)\ge n/(2d)$.
Let $x$ be an integer with
$|x-n/(2d)|\le 1/2$; say $x=n/(2d)+p$, where $-1/2\le p \le 1/2$. Choose a set $A\subseteq V_1$ with $|A|=x$,
such that all its members have at most $d$ neighbours in $V_2$. Let $B$ be the set of vertices in $V_2$
with no neighbour in $A$; then $|B|\ge n-dx$. Now
$$|A|\cdot |B|\ge x(n-dx) = (n/(2d)+p)(n-d(n/(2d)+p))=n^2/(4d)-p^2d=5t/4-p^2d.$$
But $p^2d\le t/4$ since $|p|\le 1/2$ and $d\le t$ from the hypothesis, and so $|A|\cdot |B|\ge t$.
This proves \ref{easylemma}.~\bbox

We deduce \ref{lem:onepluseps2}, slightly strengthened to the following:
\begin{thm}\label{easyEH2}
Let $H$ be a graph with $k\ge 2$ vertices. Define $\sigma = 1/(2^{k-1}-1)$.
If $G$ is an $H$-free graph with $n>1$ vertices, there is a pure pair $(A,B)$
with $|A|\cdot |B|\ge \frac{1}{45}n^{1+\sigma}$.
\end{thm}
\Proof
Since $n\ge 1$, there is a vertex either with at least $\lfloor n/2 \rfloor$ neighbours or at least
$\lfloor n/2 \rfloor$ non-neighbours,
and so we may assume that $\lfloor n/2 \rfloor< \frac{1}{45}n^{1+\sigma}$. Now $n/3\le \lfloor n/2 \rfloor$, since $n\ge 2$,
and so $1/3 < \frac{1}{45}n^{\sigma}$, that is, $n> 15^{1/\sigma}$.
Let $t= \frac{1}{45}n^{1+\sigma}$. If $t\le 5^{2^{k-2}}$ then
$$\frac{1}{45}n^{1+\sigma}\le 5^{2^{k-2}}\le 15^{1/\sigma}/3$$
and so $n^{1+\sigma}\le 15^{1/\sigma+1}$, contradicting that $n> 15^{1/\sigma}$.
Thus $t> 5^{2^{k-2}}$.

If $n\le 5kt^{e_k}+k$, then since $n/(6k)\le (n-k)/(5k)$ (because
$n> 15^{1/\sigma}$ and  $15\ge (6k)^{\sigma}$), it follows that
$n/(6k)\le t^{e_k}= 45^{-e_k}n^{(1+\sigma)e_k}=45^{-e_k}n$,
so $45^{-e_k}\ge 1/6$, a contradiction.
Thus $n> 5kt^{e_k}+k$, and so we can divide
the vertex set of $G$ into $k$ sets $V_1\ll V_k$ each of cardinality at least $5t^{e_k}$.

From \ref{easylemma} applied to the corresponding $k$-partite graph, there are sets $A,B\subseteq V(G)$,
complete or anticomplete
to each other, with $|A|\cdot |B|\ge t$, as required.~\bbox

\section{Excluding almost-bipartite graphs}\label{sec:almostbip}

We recall that a graph $H$ is \emph{almost-bipartite} if $H$ is triangle-free and there is a partition of $V(H)$ into $A, B$ 
such that $A$ is a stable set and $H[B]$ is a graph with maximum degree one. We call such a pair $(A,B)$ an \emph{almost-bipartition}.

In this section we prove \ref{mainthm}, the main theorem of the paper. It says that for every almost-bipartite graph $H$, there are $s,\vare>0$
such that for every $\vare$-bounded graph $G$, and all $0<c\le 1$, if $G$ is not $(\vare c^s, H)$-saturated then
$G$ has a $c$-sparse $(\vare c^s n, \vare n)$-pair. The difference with \ref{2csparse}
is that the latter only tells us that $G$ contains 
a $c$-sparse $(\vare c^s n, \vare c^s n)$-pair; and so far
we only know how to prove the stronger statement for almost-bipartite graphs.
Before we begin the proof (which is elaborate), it might be helpful if we sketch the main ideas.

Let us see how to do it if $H$ is actually bipartite rather than just almost-bipartite. Let $A,B$ be a bipartition, and choose
$\vare>0$ very small and $s$ very large, in terms of $H$. Now let $G$ be $\epsilon$-bounded, and let $0\le c\le 1$,
and assume $G$ has no $c$-sparse
$(\vare c^s n, \vare n)$-pair; we need to prove $G$ is $(\vare c^s, H)$-saturated. 
From \ref{lem:kcsparse}, we can arrange the constants such that
$G$ contains a $c$-sparse $(\delta c^{s_1}, |A|)$-tuple (for some constant $s_1$ much smaller than $s$). So,
take such an $a$-tuple $C_1\ll C_a$ say, where $a = |A|$. From now on we will only count copies of $H$ where for each $i$, 
the vertex
representing the $i$th vertex of $A$ is contained in $C_i$, and hope this will already give us enough copies.
This $a$-tuple is $c$-sparse, so if we pick a vertex from each at random, then with high probability the transversal we generate is stable.
This property is crucial. However, we are going to need to shrink the sets to a small fraction 
of their original size (scaled by powers of $c$) and these shrunken sets may be very dense to one another, and we might lose 
the crucial property that transversals are mostly stable. We can avoid this by choosing the original sets more carefully, 
using \ref{lem:splusone} with some 
large value of $t$,  instead of just
\ref{lem:kcsparse}; so let us do that instead. Now the edges between the 
$C_i$'s will give us no further trouble.

Most of the vertices of $G$ lie in none of $C_1\ll C_a$ (we can prove that all the $C_i$'s have cardinality at most $\vare n$);
and each vertex in $C_i$ is only adjacent to at most $\vare n$ of the outside vertices. Consequently most of the outside vertices
are only adjacent to at most $2\vare |C_i|$ vertices in $C_i$; discard the others. Actually, just discard those adjacent to more than
$2a\vare|C_i|$ vertices in $C_i$; we can afford to do this for each $i$ and still keep a good fraction of the outside vertices.

If $D$ is the set of surviving vertices outside $C_1\ll C_a$, the pair $(C_i,D)$ is not a $c$-sparse $(\vare c^s n, \vare n)$-pair, so
most vertices outside have at least $c|C_i|/(2a)$ neighbours in $C_i$; discard those that do not, for each~$i$. 
Thus $D$ still contains a constant fraction of the original vertices of $G$, where the constant depends on $H$ 
(actually, just on $a$) but not on $c$. Pick any one of those vertices, decide it is going to represent the first vertex $b_1$ 
say of $B$,
and shrink all the sets $C_i$ so that $v$ is complete to some of them and anticomplete to the others, according to 
the vertices in $A$ that $b_1$ is adjacent to in $H$. Now repeat for $b_2$, and so on; the sets $C_i$
are shrinking by factors of $c$ at each stage, but the number of choices for the next vertex in $B$ remains linear in $n$
independent of $c$. This would prove that $G$ is $(\vare c^s, H)$-saturated when $H$ is bipartite, since the shrunken sets
still have the property that random transversals are mostly stable.

How can we modify the proof to work when $H$ is almost-bipartite? Let $(A,B)$ be the almost-bipartition.
We start almost the same, applying \ref{lem:splusone} to get a $c^{s_1+t}$-sparse $(\delta c^{s_1}, 6|A|)$-tuple,
$C_1\ll C_{6a}$. (Note the 6.) Previously we filled in the vertices of $B$ one at a time, proving there were linearly 
many choices at each step. Now we fill in the edges of $H[B]$ one at a time, that is, we will add the vertices of $B$
two at a time in adjacent pairs. (We can assume that $H[B]$ is a perfect matching.) As before, we can arrange that
random transversals of the $C_i$'s are mostly stable, even after shrinking the $C_i$'s by factors of $c$; 
and that every vertex outside has a decent number of neighbours in each $C_i$ (not too many and not too few). 
So any vertex outside can play the role of any one vertex of $B$, but how do we get an edge of outside vertices to
represent an edge of $H$? If we hope an edge $uv$ of outside vertices could represent an edge $b_1b_2$
of $H$, we need $u$ to have many neighbours nonadjacent to $v$ in certain of the $C_i$'s (because $b_1$ has certain neighbours 
nonadjacent
to $b_2$ in $H$), and vice versa;
and it need not have any such neighbours. We don't know how to control things directly in this way. 

On the other hand, we can get many ($\Omega(n^2)$, regarding $c$ as a constant) edges $uv$ of outside vertices that ``disagree'' in this way
for at least one-sixth of the values of $i$. The bad news is, we can't control which one-sixth of the values this is. But there is
also good news: on that one-sixth of the values, we can shrink the sets $C_i$ to make the adjacency to $uv$
whatever we want (except, no triangles). Slightly more exactly, we can arrange that for many edges $uv$ of outside vertices, 
there exists $I\subseteq  \{1\ll 6a\}$ with $|I|= a$,
such that for each $i\in I$ there are many
vertices in $C_i$ adjacent to $u$ and not $v$, {\em and} many adjacent to $v$ and not $u$. (We won't need both sets
for a given value of $i$, because $H$ is triangle-free; but we don't yet know which set we will need.) So, pick one of 
these pairs $uv$; for the five-sixth of the values of $i$ not in $I$, we shrink $C_i$ to make it anticomplete to both $u,v$; 
and for each $i\in I$, we can
shrink the sets $C_i$ to make $C_i$ complete to $u$ and anticomplete to $v$, or vice versa, or anticomplete to both; whichever 
we want. 

Which should we choose?
We are given the power to add an edge with any adjacency we like to $I$, but with no control over the set
$I$. We can make this work as follows. If we can turn the $a$ sets corresponding to $I$ into a (blowup of) a copy of $H$
by using the edge $uv$  and shrinking the $C_i$'s for $i\in I$ appropriately (and using some of the edges added at earlier steps), 
do so;
and if not, make $I$ closer to being part of a blowup of $H$. After adding
a bounded number of edges (making sure all the added edges are anticomplete to one another, and adding each to be as useful
as possible in this way), there must be a blowup of $H$, because at every step, some $a$-subset of $\{1\ll 6a\}$ gets closer
to being in a blowup of $H$, and there are only ${6a\choose  a}$ such subsets. That is the proof; we would have shown that
there are $\Omega(n^{2k})$ induced $k$-edge matchings in $G$ (for some constant $k$), each including a submatching of 
size $b/2$ that for some $I\subseteq \{1\ll 6a\}$ with $|I|=a$, has the correct 
adjacency to a big subset of $C_i$ for each $i\in I$, and which therefore extends to $\Omega(n^a)$ copies of $H$.

Let us say this carefully. 
Let $G$ be a graph, and let $C, D \subseteq V(G)$ be disjoint. A pair $(R,S)$ of disjoint subsets of $D$ is \emph{$C$-split}
if $|N(u) \cap C|\le |N(v) \cap C|$ for each $u\in R$ and $v\in S$.

\begin{thm} \label{lem:bigsmall}
Let $G$ be a graph and $d,k \in \mathbb{N}$.
Let $C_1, \dots, C_k, D \subseteq V(G)$ with $|D|\ge 2^{k}d$.
Then there is a set
$I \subseteq \sset{1, \dots, k}$ with $|I| \geq k/2$, and
disjoint subsets $A, B$ of $D$ with $|A|, |B| = d$ such that $(A, B)$ is $C_i$-split for all $i \in I$.
\end{thm}
\Proof
We start by proving the following claim.
\\
\\
(1)  {\em Let $G$ be a graph, and $d,k\ge 1$ be integers. Let $C_1, \dots, C_k, D \subseteq V(G)$ with $|D|\ge 2^{k}d$. Then there
are disjoint subsets $A, B$ of $D$ with $|A|, |B| = d$ such that
for all $i \in \sset{1, \dots, k}$, either $(A, B)$ or $(B, A)$ is $C_i$-split.}
\\
\\
We prove this by induction on $k$. For $k = 1$, we have $|D|\ge 2d$; choose distinct $v_1\ll v_{2d}\in D$, numbered such that
$|C_1\cap N[v_j]|\ge |C_1\cap N[v_i]|$ for all  $i,j$ with $1\le i<j\le 2d$. Then
$A=\{v_1\ll v_d\}$ and $B=\{v_{d+1}\ll v_{2d}\}$
is the desired partition.

Now let $k > 1$. From the inductive hypothesis, with $d$ replaced by $2d$, it follows that there are
disjoint subsets $A, B$ of $D$ with $|A|=|B|=2d$  such that for all
$i \in \sset{1, \dots, k-1}$, either $(A, B)$ or $(B, A)$ is $C_i$-split. Choose $j\in \mathbb{N}$ maximum such that
\begin{itemize}
\item the number of vertices in $A$ with at most $j-1$ neighbours in $C_k$ is at most $d$, and
\item the number of vertices in $B$ with at most $j-1$ neighbours in $C_k$ is at most $d$.
\end{itemize}
(This is possible, because for $j=0$ both bullets are true, and for $j\ge 2d$ both are false, since $d\ge 1$.) By exchanging
$A,B$ if necessary, we may assume that the number of vertices in $A$ with at most $j$ neighbours in $C_k$ is more than $d$.
Consequently there is a subset $A'\subseteq A$ with $|A'|=d$ such that each vertex in $A'$ has at most $j$ neighbours
in $C_k$; and a subset $B'\subseteq B$ with $|B'|=d$ such that each vertex in $B'$ has at least $j$ neighbours in $C_k$.
Thus $(A',B')$ is $C_k$-split, and for $1\le i\le k-1$, either $(A', B')$ or $(B', A')$ is $C_i$-split.
This proves (1).

\bigskip

Now we apply (1). We may assume that $d,k\ge 1$; we apply (1) to $C_1, \dots, C_k$ and obtain $A, B$ as in (1).
Then one of $(A,B), (B,A)$ satisfy the theorem.
This proves \ref{lem:bigsmall}.~\bbox

\begin{thm}\label{cbound}
Let $0<\epsilon\le 1/4$, let $s \in \mathbb{N}$ with $k \geq 1$, let $0\le c\le 1$, and let
$G$ be an $\vare$-bounded graph with $n>0$ vertices, and with no $c$-sparse $(\vare c^s n, \vare n)$-pair.
Then $n>\epsilon^{-1}$, and $c < 2\epsilon$.
\end{thm}
\Proof Let $v$ be a vertex of $G$. Then $1 \leq |N[v]| < \epsilon n$ since $G$ is $\vare$-bounded, and the
first statement follows.

For the second statement,  since $\vare \leq 1/4$, it follows that $n > 4$.  Let $A \subseteq V(G)$ with
$|A| = \floor{n/2}$. Since $|A| \geq n/4\ge \vare c^s n$, it follows that $(A, V(G) \setminus A)$ is not a $c$-sparse pair,
and so $|E(A, V(G) \setminus A)| > c |A|\cdot|V(G) \setminus A|$. Therefore, there is a vertex $v$ in $A$ with
at least $c |V(G) \setminus A|$ neighbours in $V(G) \setminus A$, so $v$ has degree at least $c n/2$.
Since every vertex has degree less than $\vare n$, it follows that $cn/2 < \vare n$,
and the second statement follows. This proves \ref{cbound}.~\bbox

This is used to prove the main lemma, that we can get many edges that disagree on one-sixth of the values of $i$.

\begin{thm}\label{onethird}
Let $k,s\ge 1$ be integers, let $0<\vare\le k^{-1} 2^{-k-5}$, and let 
$G$ be an $\vare$-bounded graph with $n$ vertices.
Let $0< c\le 1$  such that $G$ does not have a $c$-sparse $(\vare c^s n, \vare n)$-pair.
Let $C_1\ll C_k, D$ be pairwise disjoint 
subsets of $V(G)$  such that $|C_i|\ge 2\vare c^s n$ 
for $1\le i\le k$, and $|D|\ge 3n/4$. 
Let $E^*$ be the set of edges $uv$ of $G[D]$ such that for at least $k/6$ values of $i\in \{1\ll k\}$,
there are at least $c^2 |C_i|$ vertices in $C_i$ adjacent to $u$ and not to $v$, and at least $c^2 |C_i|$ adjacent to $v$ 
and not to $u$.
Then $|E^*|\ge c^22^{-2k-9}n^2$.
\end{thm}
\Proof
For $1\le i\le k$, let $D_i$ be
the set of vertices in $D$ with at most $c|C_i|$ neighbours in $C_i$; and
let $F_i$ be
the set with at least $8k \vare |C_i|$ neighbours in $C_i$.
Then $|D_i| \leq \vare n\le n/(8k)$, since $(C_i, D_i)$ is not a $c$-sparse
$(\vare c^s n, \vare n)$-pair; and $|F_i| \leq n/(8k)$,
since
$$8k\vare |C_i|\cdot|F_i| \leq |E(C_i, D)| \leq \vare |C_i|n.$$
It follows that $|D_i \cup F_i| \leq n/(4k)$ for $1\le i\le k$, and
so 
$$D' = D \setminus \bigcup_{i \in \sset{1, \dots, k}}(D_i \cup F_i)$$ 
satisfies
$|D'| \geq |D|-n/4\ge n/2$.
Let $d=\floor{2^{-k-2}n}$. Since $2^{-k-2}n\ge\vare n> 1$ by \ref{cbound},  it follows that $d\ge 2^{-k-3}n$.

Since $|D'|\ge 2^{k+1}d$, \ref{lem:bigsmall} implies there exist disjoint
subsets $A,B\subseteq D'$ with $|A|=|B|=2d$, and $I\subseteq \{1\ll k\}$ with $|I|\ge k/2$  such that 
$(B, A)$ is $C_i$-split for all $i \in I$.

For $uv \in E(G)$ with $u, v \in D'$, and $i\in I$,
we write $u \rightarrow_i v$ if $|(N(v) \setminus N(u)) \cap C_i| < c^2 |C_i|$ (note that possibly $u \rightarrow_i v$ and
$v \rightarrow_i u$ both hold). If neither of $u \rightarrow_i v$,
$v \rightarrow_i u$ hold, we say $uv$ is {\em $i$-incomparable}.
If there are at least $|I|/3$ values of $i\in I$ such that $u \rightarrow_i v$, we write
$u \rightarrow v$. (Again, possibly $u \rightarrow v$ and 
$v \rightarrow u$ both hold.) 
Thus, any edge $uv$ that is $i$-incomparable for at least $|I|/3$ values of $i\in I$ belongs to $E^*$,
and in particular, any edge $uv$ for which $u \not\rightarrow v$ and $v \not\rightarrow u$ belongs to $E^*$.
\\
\\
(1) {\em We may assume that there is a vertex $u \in A$ such that the set $U = \sset{v \in A : u \rightarrow v}$
satisfies $|U| \geq cd/4$.}
\\
\\
Choose disjoint subsets $A_1, A_2$ of $A$, both of cardinality $d$. 
Since $d\ge \vare n$, and $G$ has no $c$-sparse $(\vare c^s n, \vare n)$-pair, it follows that the pair $(A_1, A_2)$ is not
$c$-sparse. Hence $|E(A_1,A_2)|\ge cd^2$.  For each edge $uv$ with $u\in A_1$ and $v\in A_2$, either
$u \rightarrow v$, or $v \rightarrow u$, or $uv\in E^*$; and if $|E^*|\ge c^22^{-2k-9}n^2$ we are done. So we may assume 
(exchanging $A_1, A_2$
if necessary) that there are at least $(cd^2-c^22^{-2k-9}n^2)/2\ge cd^2/4$ edges $uv$ with $u\in A_1$ and $v\in A_2$, 
such that $u \rightarrow v$.
Hence there exists $u\in A_1$ such that $u \rightarrow v$ for at least $cd/4$ values of $v\in A_2$. 
This proves (1).
\bigskip

Choose $u$ and $U$ as in (1).
For $i\in I$, since $u\notin F_i$ and $i\in I$, it follows that
$|N(u) \cap C_i| \leq 8k\vare |C_i|$; and consequently
$$|C_i \setminus N(u)| \geq (1-8k\vare) |C_i|\ge  \vare c^s n.$$

For $i \in I$, let $B_i$ be the set of vertices in $B$ with at most $c|C_i \setminus N(u)|$ neighbours in $C_i \setminus N(u)$.
Since $(C_i \setminus N(u), B_i)$ is a $c$-sparse pair,
and $|C_i \setminus N(u)|\ge \vare c^s n$, it follows that 
$|B_i| \leq \vare n$, for all $i \in I$.
Let $B' = B \setminus \left(\bigcup_{i \in I} B_i\right)$; then
$|B'| \geq |B| - k \vare n \ge d/2$.
Since $|U|\ge cd/4\ge \vare c^s n$ and $|B'|\ge \vare n$, it follows that
$U, B'$ is not $c$-sparse, and so
$$|E(U,B')| \geq c|U|\cdot|B'| \geq c^2d^2/8\ge c^22^{-2k-9}n^2.$$
We claim that $E(U,B')\subseteq E^*$, and the result will follow. 
\\
\\
(2) {\em $E(U,B')\subseteq E^*$.}
\\
\\
Let $vw$ be an edge with $v\in U$ and $w\in B'$; and let $I' = \sset{i \in I: u \rightarrow_i v}$.
We need to show that $vw\in E^*$; and to prove this, it suffices to show that for each $i\in I'$, 
$v \not\rightarrow_i w$ and $w \not\rightarrow_i v$, and so $v,w$ are $i$-incomparable.
Thus, let $i\in I'$. Now $u$ has at most $8k \vare |C_i|\le |C_i|/2$ neighbours in $C_i$, since $u\notin F_i$; 
and so $|C_i \setminus N(u)|\ge |C_i|/2$.
But $w$ has at least $c|C_i \setminus N(u)|\ge c|C_i|/2$ neighbours in
$C_i \setminus N(u)$ since $w \in B'$, and $v$ has at most $c^2|C_i|$ neighbours in $C_i \setminus N(u)$ since $u\rightarrow_i v$.
Hence there are at least $c|C_i|/2-c^2|C_i|\ge c^2|C_i|$ vertices in $C_i$ that are adjacent to $w$ and not to $v$,
since $c<2\vare\le 1/4$ by \ref{cbound}, and
it follows that $v \not\rightarrow_i w$. 
Moreover, we recall that
$(B, A)$ is $C_i$-split since
$i\in I'\subseteq I$; and since
$v \in A, w \in B$, it follows that $|N(w) \cap C_i| \leq |N(v) \cap C_i|$. Consequently
$$|(N(v) \setminus N(w)) \cap C_i| \geq |(N(w) \setminus N(v)) \cap C_i| \geq c^2 |C_i|,$$
and so $w \not\rightarrow_i v$. This proves (2).

\bigskip

From (2), $|E^*|\ge |E(U,B')| \geq c^22^{-2k-9}n^2.$ This proves \ref{onethird}.~\bbox

An {\em induced matching} in a graph $G$ is a subset $M\subseteq E(G)$  such that 
for all distinct $e,f\in M$, $e,f$ have no common end and 
both ends of $e$ are nonadjacent to both ends of $f$. We write $V(M)$ to denote 
the set of ends of members of $M$.
A {\em blockade} in a graph $G$ is a set $\mathcal{C}=\{C_1\ll C_k\}$ of pairwise disjoint 
subsets of $V(G)$. (We used the same term to mean something slightly different in \cite{trees}.)
We write $V(\mathcal{C})=C_1\cup\cdots\cup C_k$.
If $C_i'\subseteq C_i$ for $1\le i\le k$, 
we call $\mathcal{C}'=\{C_1'\ll C_k'\}$
a {\em contraction} of $\mathcal{C}$; and for $\delta> 0$, if $|C_i'|\ge \delta|C_i|$ for $1\le i\le k$.
we call $\mathcal{C}'$ a {\em $\delta$-contraction} of $\mathcal{C}$.

Now let 
$M$ be an induced matching in $G\setminus V(\mathcal{C})$. We say $\mathcal{C}$ is
{\em $M$-pure} if 
\begin{itemize}
\item for $1\le i\le k$, and for each $v\in V(M)$, $v$ is either complete or anticomplete to $C_i$; and
\item for each $e=uv\in M$ and $1\le i\le k$, not both $u,v$ are complete to $C_i$.
\end{itemize}
Let $\mathcal{C}$ be $M$-pure.
For each $e=uv\in M$,
let $P,Q,R$ be respectively the sets of $i\in \{1\ll k\}$ such that $u$ is complete to $C_i$, $v$ is complete to $C_i$, and neither;
then $(P,Q,R)$ is a partition of $\{1\ll k\}$, and we call $(P,Q,R)$ and $(Q,P,R)$ the {\em supports} of $e=uv$.
The number of distinct supports of edges in $M$ is called the {\em richness} of $M$ on $\mathcal{C}$.
(More precisely, the richness is the number of partitions $(P,Q,R)$ of $M$ such that $(P,Q,R)$ is a support of an edge of $M$.)

If $\mathcal{C}=\{C_1\ll C_k\}$ is a blockade and $I\subseteq \{1\ll k\}$, then $\{C_i:i\in I\}$ is called a 
{\em sub-blockade.} If $\mathcal{C}$ is $M$-pure, then so are its sub-blockades. We say $M$ is {\em complete} on $\mathcal{C}$
if every partition of $\{1\ll k\}$ into three parts is the support of an edge of $M$.

Let $a\in \mathbb{N}$, and let $k=6a$; and let $\mathcal{C}=\{C_1\ll C_k\}$ be a blockade. Then the {\em worth} of $M$ on $\mathcal{C}'$ is $\infty$ if $M$ is complete on some sub-blockade of cardinality $a$; and otherwise
the worth of $M$ is the sum, over all sub-blockades $\mathcal{C}'$ of cardinality $a$, 
of the richness of $M$ on $\mathcal{C}'$. Thus, the worth of $M$ is either less than $3^a{6a\choose a}$ or $\infty$.

Let us say an induced matching $M$ in $G\setminus V(\mathcal{C})$ is {\em $\mathcal{C}$-successful}  
if there is an $M$-pure
$c^{2|M|}$-contraction $\mathcal{C}'$
of $\mathcal{C}$  such that $M$ has worth at least $|M|$ over $\mathcal{C}'$.

\begin{thm}\label{growmatching}
Let $a\ge 1$ be an integer, let $k=6a$, and let $R=3^a{k\choose a}$; let $s\ge 2R$ be an integer, and let $0<\vare\le 6^{-k}$.
Let $G$ be an $\vare$-bounded graph with $n$ vertices.
Let $0<c\le 1$  such that $G$ does not have a $c$-sparse $(\vare c^s n, \vare n)$-pair.
Let $\mathcal{C}=\{C_1\ll C_k\}$ be a blockade in $G$,
such that $2\vare c^{s-2R} n\le |C_i|\le \vare n$
for $1\le i\le k$. Then for $0\le m\le R$ there are at least $2^{-(2k+9)m}c^{2m}n^{2m}/m!$
$\mathcal{C}$-successful induced matchings in $G\setminus V(\mathcal{C})$ of cardinality $m$.
\end{thm}
\Proof
We show first:
\\
\\
(1) {\em Let $0<m\le 3^a{6a\choose a}$, and let $M$ be a $\mathcal{C}$-successful induced matching in $G\setminus V(\mathcal{C})$ of cardinality $m-1$.
Then there are at least $c^22^{-2k-9}n^2$ $\mathcal{C}$-successful induced matchings in $G\setminus V(\mathcal{C})$ of cardinality $m$
that include $M$. }
\\
\\
Since $M$ is $\mathcal{C}$-successful, there is an $M$-pure 
$c^{2|M|}$-contraction $\mathcal{C}'=(C_1'\ll C_k')$
of $\mathcal{C}$  such that $M$ has worth at least $|M|$ on $\mathcal{C}'$. Let $D$ be the set of all vertices in 
$V(G)\setminus V(\mathcal{C})$ that are anticomplete to $V(M)$. Thus $|D|\ge n-(k+2m)\vare n$, since each 
$|C_i|\le \vare n$ and each vertex in $V(M)$ has at most $\vare n$ neighbours; and so $|D|\ge 3n/4$, since
$(k+2m)\vare\le 1/4$ because $m\le 3^a{6a\choose a}$ and $\vare \le 6^{-k}$.
Now each $|C_i'|\ge c^{2m-2}|C_i|\ge 2\vare c^{s} n$, since $m\le R$ and $s\ge 2R$.
Since $(k+2m)\vare n\le n/4$, and $\vare\le 6^{-k}\le k^{-1} 2^{-k-5}$,
we can apply \ref{onethird} to $\mathcal{C}'$ and $D$. With $E^*$ defined as in \ref{onethird}, 
we deduce that $|E^*|\ge c^22^{-2k-9}n^2$.
Let $e=uv\in E^*$. We claim that $M\cup \{e\}$ is a $\mathcal{C}$-successful matching.

From the definition of $E^*$, there exists $I\subseteq \{1\ll k\}$ with $|I|=a$ such that for each $i\in I$,
there are at least $c^2 |C_i|$ vertices in $C_i$ adjacent to $u$ and not to $v$, and at least $c^2 |C_i|$ adjacent to $v$
and not to $u$. There are two cases:
\begin{itemize}
\item If $M$ has worth $\infty$ on $\mathcal{C}'$, then for each $i\in \{1\ll k\}$,
let $C_i''$ be the set of vertices in $C_i'$ nonadjacent to both $u,v$, and $\mathcal{C}''=\{C_1''\ll C_k''\}$. Then
$|C_i''|\ge |C_i'|-2\vare n\ge c^2 |C_i'|$ for each $i$, and $\mathcal{C}''$ is $(M\cup \{e\})$-pure;
and $M\cup\{e\}$ has worth $\infty$ on $\mathcal{C}''=\{C_1''\ll C_k''\}$.
\item If $M$ has finite worth on $\mathcal{C}'$, then in particular $M$ is not complete on the sub-blockade $\{C'_i:i\in I\}$
of $\mathcal{C}'$. Choose a partition $(P,Q,R)$ of $I$ that is not a support of any edge
in $M$; for each $i\in P$, let $C_i''$ be the set of vertices in $C_i'$ adjacent to $u$ and not to $v$; for each $i\in Q$,
let $C_i''$ be the set of vertices in $C_i'$ adjacent to $v$ and not to $u$; and for each $i\in \{1\ll k\}\setminus (P\cup Q)$,
let $C_i''$ be the set of vertices in $C_i'$ nonadjacent to both $u,v$. Then again, for each $i$, $|C_i''|\ge c^2 |C_i'|$,
and $\mathcal{C}''$ is $(M\cup \{e\})$-pure.
Moreover, the worth of $M\cup\{e\}$ on $\mathcal{C}''=\{C_1''\ll C_k''\}$ is strictly more than the worth of $M$ on $\mathcal{C}'$
(because
the richness on the sub-blockade defined by $I$ increased).
\end{itemize}
In both cases,
$\mathcal{C}''$ is a $c^2$-contraction of $\mathcal{C}'$ and hence a $c^{2m}$-contraction of $\mathcal{C}$; and the worth
of $M\cup \{e\}$ on $\mathcal{C}''$ is either $\infty$, or at least one more than the worth of $M$ on $\mathcal{C}'$ 
and so in either case
$M\cup \{e\}$ has worth at least $m$ on $\mathcal{C}''$. Consequently $M\cup \{e\}$ is a $\mathcal{C}$-successful matching,
for at least $c^22^{-2k-9}n^2$ edges $e$. This proves (1).

\bigskip

For $m\ge 0$, let $f_m$ denote the number of $\mathcal{C}$-successful induced matchings of cardinality $m$.
We must show that $f_m\ge 2^{-(2k+9)m}c^{2m}n^{2m}/m!$.
We proceed by induction on $m$; the claim is true if $m=0$, so we assume $m>0$ and the claim holds for $m-1$.
Since every $m$-edge matching includes only $m$ matchings of cardinality $m-1$, it follows from (1) that
$$mf_m\ge c^22^{-2k-9}n^2 f_{m-1}\ge c^22^{-2k-9}n^22^{-(2k+9)(m-1)}c^{2m-2}n^{2m-2}/(m-1)!,$$ 
and so $f_m\ge 2^{-(2k+9)m}c^{2m}n^{2m}/m!$. This proves \ref{growmatching}.~\bbox

We also need a lemma about sparse $k$-tuples.

\begin{thm}\label{trimtuple2}
Let $G$ be a graph, let $c>0$, and let $k\ge 2$ be an integer. Let $C_1\ll C_k$ be pairwise disjoint subsets of $V(G)$,
pairwise $c$-sparse.
Then there exist $B_i\subseteq C_i$
for $1\le i\le k$ such that $|B_i|\ge |C_i|/2$ for $1\le i\le k$, and for all distinct $i,j\in k$, every vertex in $B_i$
has at most $4kc|B_j|$  neighbours in $B_j$.
\end{thm}
\Proof
For each $j\ne i$, let $A_{i,j}$ be the set of vertices in $C_i$ with at least $2kc|C_j|$ neighbours in $C_j$. Since there are only
$c|C_i|\cdot|C_j|$ edges between $C_i$ and $C_j$, it follows that $|A_{i,j}|\le c|C_i|\cdot|C_j|/(2kc|C_j|) = |C_i|/(2k)$, and
so the union of all the sets $A_{i,j} (j\ne i)$ has cardinality at most $|C_i|/2$. Let $B_i$ be its complement in $C_i$.
Thus $|B_i|\ge |C_i|/2$, and for all distinct $i,j$, each vertex in $B_i$ has at most $2kc|C_j|\le 4kc|B_j|$ neighbours in $C_j$,
and therefore has at most that many in $B_j$. This proves \ref{trimtuple2}.~\bbox

Now we are ready to prove our main theorem:

\begin{thm}\label{almostbip}
Let $H$ be almost-bipartite. Then there exist $s \in \mathbb{N}$ and $\vare > 0$ such that for all $c$ with $0<c\le 1$, if $G$ 
is $\vare$-bounded, then either $G$ is $(\vare c^s, H)$-saturated, or $G$ has a $c$-sparse $(\vare c^s |G|, \vare|G|)$-pair. 
\end{thm}
\Proof
There is a ``universal'' almost-bipartite graph, defined as follows. Let $a>0$ be an integer; and let $A$ be a set of cardinality $a$.
For each partition $(P,Q,R)$ of $A$, where either $P=\emptyset$ or the least member of $P$ is less than all members of $Q$,
take an edge $uv$ of new vertices, and make $u$ adjacent to all $p\in P$ and $v$ adjacent to all $q\in Q$. We call the result
$H_a$; it is almost-bipartite, and it is easy to see that every almost-bipartite graph is an induced subgraph of $H_a$ for some $a$.
(The restriction on $(P,Q,R)$ is to ensure that no two of the added edges have a common support.) Thus, to prove \ref{almostbip}
in general, it suffices to prove it when $H=H_a$, for each $a$.

Let $k=6a$, $R= 3^a{k\choose a}$, and $t=2R+1$. Let \ref{lem:splusone} hold with $k,t,H,S, \vare$ replaced by $k,t, H_a, S, \vare'$ respectively.
Define $p=2^{-(2k+9)R}/R!$, and $q=(\vare'/2)^{k}/2$. 
Let $s=2R+k(2R+S)$ and 
$$\vare = \min\left(pq,\vare'/2,1/(8k^2a^2), 6^{-k}\right).$$ 
We claim that $s,\vare$ satisfy the theorem. For let $0<c\le 1$, and let 
$G$
be $\vare$-bounded, with $n>0$ vertices. 

We may assume that
$G$ has no $c$-sparse $(\vare c^s n, \vare n)$-pair; and now we must show that $G$ is $(\vare c^s, H)$-saturated.
From \ref{cbound}, $c<2\vare$.
We may assume that $G$ is not $\left(\vare' c^{S}, H\right)$-saturated, since $S\le s$ and $\vare'\ge \vare$; 
so from \ref{lem:kcsparse}, $G$
contains a $c^{s'+t}$-sparse $(\vare' c^{s'}, k)$-tuple $C_1\ll C_{k}$ for some $s'\le S$.

By \ref{trimtuple2}, with $c$ replaced by $c^{s'+t}$, there exists $B_i\subseteq C_i$ with $|B_i|\ge |C_i|/2$ for $1\le i\le k$, 
such that
for all distinct $i,j\in \{1\ll k\}$, every vertex in $B_i$
has at most $4kc^{s'+t}|B_j|$  neighbours in $B_j$. Let $\mathcal{B}=(B_1\ll B_k)$.
For all distinct $i,j\in \{1\ll k\}$, since $G$ has no $c$-sparse $(\vare c^s n, \vare n)$-pair, and the pair 
$(B_1, B_2)$ is $c$-sparse, and $|B_i|\ge \vare' c^{s'}n/2\ge \vare c^s n$ 
(since $\vare\le \vare'/2$ and $s'\le S\le s$), it follows that $|B_j|\le \vare n$; 
and so $|B_1|\ll |B_{k}|\le \vare n$. 
By \ref{growmatching}, there are at least $pc^{2R}n^{2R}$ $\mathcal{B}$-successful
induced matchings in $G\setminus V(\mathcal{B})$ of cardinality $R$.
\\
\\
(1) {\em For every $\mathcal{B}$-successful induced matching $M$ in $G\setminus V(\mathcal{C})$ of cardinality $R$,
there are at least $qc^{k(2R+S)}n^{k}$ choices of $(x_1\ll x_k)$  such that $x_i\in C_i$ for $1\le i\le k$, and
the subgraph induced on 
$$V(M)\cup \{x_1\ll x_{k}\}$$ 
contains $H_a$.}
\\
\\
Let $M$ be a $\mathcal{B}$-successful induced matching of cardinality $R$. Hence there is an $M$-pure
$c^{2R}$-contraction $\mathcal{A} = (A_1\ll A_{k})$
of $\mathcal{B}$  such that $M$ has worth at least $|M|$ on $\mathcal{A}$. Since for all distinct $i,j\in \{1\ll k\}$,
every vertex in $B_i$ has at most $4kc^{s'+t}|B_j|$ neighbours in $B_j$, there are at most 
$$4kc^{s'+t}|A_i|\cdot|B_j|\le 4kc^{s'+t-2R}|A_i|\cdot|A_j|\le 4kc|A_i|\cdot|A_j|$$
edges between $A_i$ and $A_j$, since $s'+t-2R\ge 1$; that is, $A_1\ll A_k$ is $4kc$-sparse.

Now since $|M|=R$, it follows that $M$ has worth $\infty$
on $\mathcal{A}$,
and so there exists $I\subseteq \{1\ll k\}$ with $|I|=a$ such that $M$ is complete on the sub-blockade $(A_i:i\in I)$. 
Let $N$ be the product of the cardinalities of $A_1\ll A_k$.
There are $N$ choices of a sequence $(x_1\ll x_k)$ such that $x_i\in A_i$ for $1\le i\le k$. For all distinct $i,j\in I$,
there are only $4kcN$ such choices in which $x_i, x_j$ are adjacent, since $(A_i, A_j)$ is $4kc$-sparse.
Hence there are at least 
$(1-4ka^2c)N$ 
choices such that $\{x_i:i\in I\}$ is stable; and since $1-4ka^2c\ge 1/2$ (because $c<2\vare\le 1/(8k^2a^2)$),
this number is at least $N/2$.
For each 
such choice of $(x_1\ll x_k)$, the subgraph induced on $V(M)\cup \{x_i:i\in I\}$
contains $H_a$, since $M$ is complete on the sub-blockade $(A_i:i\in I)$. Since for each $i$, 
$|A_i|\ge c^{2R}|B_i|\ge c^{2R+s'}\vare' n/2$,
and $s'\le S$, it follows that
$$N\ge c^{k(2R+S)}(\vare'/2)^kn^k=2q c^{k(2R+S)}n^k,$$
and this proves (1).

\bigskip

Multiplying the number of choices for $M$ and the number of choices for $(x_1\ll x_{k})$ (for each $M$), we deduce that
altogether there are at least $pqc^{2R+k(2R+s')r}n^{2R+k}$
distinct induced subgraphs of $G$, each with $k+2R$ vertices, and each containing $H_a$.
Since $pq\ge \vare$ and $2R+k(2R+s')\le s$, it follows that
there are at least $\vare c^s n^{2R+k}$ such subgraphs.
But
each induced subgraph of $G$ isomorphic to $H_a$ is contained in at most $n^{k+2R-|H_a|}$ induced subgraphs of $G$ with $k+2R$ 
vertices, and so
there are at least 
$\vare c^s n^{|H_a|}$
distinct copies of $H_a$ in $G$. This proves \ref{almostbip}.~\bbox


{\em Subdividing} an edge $uv$ of a graph $H'$ means replacing $uv$ by a path, whose internal vertices are not in $V(H')$
and have degree two in the new graph.
A graph $H$ is a {\em $(\ge 1)$-subdivision} of a graph $H'$ if $H$ arises from $H'$ by subdividing each edge at least once; that
is, replacing each edge (one at a time) by a path of length at least two. 
\ref{almostbip} implies the following. 

\begin{thm}\label{subdivision}
Let $H'$ be a graph, and let $H$ be a $(\ge 1)$-subdivision of $H'$. Then there exist $s \in \mathbb{N}$ and $\vare > 0$ such that 
for all $c > 0$, if $G$ is $\vare$-bounded and $H$-free, then $G$ has a $c$-sparse $(\vare c^s |G|, \vare|G|)$-pair. 
\end{thm}
\Proof
By  \ref{almostbip}, it suffices to prove that $H$ is almost-bipartite. We first note that since $H$ is a $(\ge 1)$-subdivision
of some graph, it follows that $H$ arises from some graph $H'$ by subdividing every edge of $H'$ either once or twice.

Now $V(H')$ is a stable set in $H$, and every vertex of $H \setminus V(H')$ has  
degree at most one in $H \setminus V(H')$ (since it is part of a path of length two or three with ends in $V(H')$). 
Since $H$ is a $(\ge 1)$-subdivision of $H'$, it follows that $H$ is triangle-free, and therefore $H$ is almost-bipartite. 
This proves \ref{subdivision}.~\bbox


\begin{thebibliography}{99}
\bibitem{ehfive} M. Chudnovsky, J. Fox, A. Scott, P. Seymour and S. Spirkl. ``Towards Erd\H{o}s-Hajnal for graphs with no 5-hole'', 
{\em Combinatorica} 39 (2019), 983--991, {\tt arXiv:1803.03588}.
\bibitem{trees} M. Chudnovsky, A. Scott, P. Seymour and S. Spirkl, ``Pure pairs. I. Trees and linear anticomplete pairs'',
submitted for publication, {\tt arXiv:1809.00919}.
\bibitem{sparse} M. Chudnovsky, A. Scott, P. Seymour and S. Spirkl, ``Pure pairs. II. Excluding all subdivisions of a graph'',
submitted for publication, {\tt arXiv 1804.01060}.
\bibitem{fox} D. Conlon, J. Fox and B. Sudakov, ``Recent developments in graph Ramsey theory'', 
in: {\em Surveys in Combinatorics 2015}, Cambridge University Press, 2015, 49--118.
\bibitem{EH0} P. Erd\H{o}s and A. Hajnal, ``On spanned subgraphs of graphs'',
{\em Graphentheorie und Ihre Anwendungen} (Oberhof, 1977), \verb++{www.renyi.hu/\raisebox{-1ex}{\textasciitilde}p\_erdos/1977-19.pdf}.
\bibitem{EH}  P. Erd\H{o}s and A. Hajnal, ``Ramsey-type theorems'',
{\em  Discrete Applied Math.} {\bf 25} (1989), 37--52.
\bibitem{EHP} P. Erd\H{o}s, A. Hajnal and J. Pach, ``A Ramsey-type theorem for bipartite graphs'', {\em Geombinatorics}
  10 (2000), 64--68.
\bibitem{foxsudakov} J. Fox and B. Sudakov, ``Induced Ramsey-type theorems'', \emph{Advances in Math.}, {\bf 219} (2008), 1771--1800.
\bibitem{cats} A. Liebenau, M. Pilipczuk, P. Seymour and S. Spirkl, ``Caterpillars in Erd\H{o}s-Hajnal'',
{\em J. Combinatorial Theory, Ser. B}, 136 (2019), 33--43, {\tt arXiv:1810.00811}.
\bibitem{rodl} V. R\"odl, ``On universality of graphs with uniformly distributed edges'',
{\em Discrete Math.} {\bf 59} (1986), 125--134.
\end{thebibliography}
\end{document}